\documentclass[11pt,oneside,english]{amsart}

\usepackage{amssymb}
\usepackage{mathtools}
\usepackage{empheq}[2004/10/10]
\usepackage[T1]{fontenc}
\usepackage[latin9]{inputenc}
\usepackage{babel}
\usepackage{graphicx}
\usepackage{xspace}
\usepackage{array,hhline}

\usepackage[unicode=true, pdfusetitle,
 bookmarks=true,bookmarksnumbered=false,bookmarksopen=false,
 breaklinks=false,pdfborder={0 0 0},backref=false,colorlinks=false]
 {hyperref}
\hypersetup{%
pdfsubject={MSC 2000 Primary: 13P10; 90C57; Secondary: 65H10;12Y05; 90C27; 68R05},
pdfkeywords={order ideal polytope, border bases, Gröbner bases, combinatorial
optimization},
pdfauthor={Gábor Braun and Sebastian Pokutta}}

\numberwithin{equation}{section} 
\numberwithin{figure}{section} 
\theoremstyle{plain}
\newtheorem{thm}{Theorem}[section]
\newtheorem{corollary}[thm]{Corollary}
  \theoremstyle{definition}
  \newtheorem{defn}[thm]{Definition}
 \theoremstyle{definition}
  \newtheorem{example}[thm]{Example}
  \theoremstyle{remark}
  \newtheorem{rem}[thm]{Remark}
  \theoremstyle{plain}
  \newtheorem{lem}[thm]{Lemma}
  \theoremstyle{plain}
  \newtheorem{inneralgorithm}[thm]{Algorithm}
  \theoremstyle{plain}
  \newtheorem{prop}[thm]{Proposition}

\usepackage[bitstream-charter]{mathdesign}
\usepackage{fullpage}

\newenvironment{algorithm}[4]{%
  \begin{inneralgorithm}[{#1\textemdash#2}]\mbox{}%
    \let\saveditem\item
    \renewcommand{\item}{%
      \let\item\saveditem
      \begin{description}
      \item[Input] #3
      \item[Output] #4
      \end{description}
    \begin{enumerate}\item}}%
    {%
    \end{enumerate}%
  \end{inneralgorithm}}

\newtheoremstyle{problem@}
  {}
  {}
  {}
  {}
  {\scshape}
  {:}
  { }
  {\thmnote{#3}}
\theoremstyle{problem@}
\newcounter{problem}
\newtheorem{problem@}{}[problem]
\makeatletter
\newenvironment*{problem}[1]{%
  \def\@currentlabelname{\textsc{#1}}
  \begin{problem@}[{#1}]}{\end{problem@}}
\makeatother

\newcommand*{\rk}[1]{\operatorname{rk}(#1)}

\newcommand*{\conv}[1]{\operatorname{conv}(#1)}
\newcommand{\cube}[1]{{[0,1]}\sp{#1}}
\newcommand{\Z}{\mathbb{Z}}
\newcommand{\Q}{\mathbb{Q}}

\newcommand{\Pclass}{\textrm{P}\xspace}
\newcommand{\NP}{\textrm{NP}\xspace}
\newcommand{\N}{\mathbb{N}}

\newcommand{\T}{\mathbb{T}}
\DeclareMathOperator{\supp}{supp}
\newcommand{\oi}{\mathcal{O}}
\DeclareMathOperator{\LT}{LT}
\DeclareMathOperator{\LC}{LC}
\DeclareMathOperator{\LF}{LF}
\newcommand*{\textprogram}[1]{\textnormal{\texttt{#1}}}
\newcommand*{\textalgorithm}[1]{\textnormal{\texttt{#1}}}
\newcommand{\bbasis}{\ensuremath{\operatorname{\textalgorithm{BBasis}}}}

\newcommand{\lstabspan}{\ensuremath{\operatorname{\textalgorithm{LStabSpan}}}}
\newcommand{\basisT}{\ensuremath{\operatorname{\textalgorithm{BasisTransformation}}}}

\newcommand{\finalRed}{\ensuremath{\operatorname{\textalgorithm{FinalRed}}}}
\newcommand{\head}{\ensuremath{\operatorname{\textalgorithm{Head}}}}
\newcommand{\gausEl}{\ensuremath{\operatorname{\textalgorithm{GaussEl}}}}
\newcommand{\pring}{K[\mathbb{X}]}
\newcommand{\vgen}[1]{{\left\langle #1\right\rangle}\sb{K}}
\newcommand{\igen}[1]{{\left\langle #1\right\rangle}\sb{\pring}}

\begin{document}
\renewcommand{\sectionautorefname}{Section}
\renewcommand{\subsectionautorefname}{Subsection}
\renewcommand{\equationautorefname}{Condition}

\title{A polyhedral approach to computing border bases}

\author{Gábor Braun}
\address{Alfréd Rényi Institute of Mathematics\\
  Hungarian Academy of Sciences\\
  Reáltanoda u 13\textendash15\\
  1053\\
  Hungary}
\email{braung@renyi.hu}
\thanks{Research partially supported by Hungarian Scientific Research Fund,
  grant No. K 67928.}

\author{Sebastian Pokutta}
\thanks{Research partially supported by German Research Foundation
  (DFG) funded SFB 805.}

\address{Fachbereich Mathematik \\ Technische Universität Darmstadt \\ Germany}

\email{pokutta@mathematik.tu-darmstadt.de}


\date{\today}

\begin{abstract}
  Border bases can be considered to be the natural extension of Gröbner bases
  that have several advantages. Unfortunately, to date the classical border
  basis algorithm relies on (degree-compatible) term orderings and implicitly
  on reduced Gröbner bases.
  We adapt the classical border basis
  algorithm to allow for calculating border bases for arbitrary
  degree-compatible order ideals,
  which is \emph{independent} from term orderings.
  Moreover,
  the algorithm also supports calculating degree-compatible order ideals
  with \emph{preference} on contained elements, even though finding a
  preferred order ideal is NP-hard.
  Effectively we retain degree-compatibility only to
  successively extend our computation degree-by-degree.
  
  The adaptation is based on our polyhedral characterization:
  order ideals that support a border basis correspond one-to-one to
  integral points of the order ideal polytope.
  This establishes a crucial connection between the ideal and the combinatorial structure of the associated factor spaces.
\end{abstract}
\maketitle

\section{Introduction}

Gröbner bases are fundamental tools in commutative algebra to model
and perform important operations on ideals such as intersection, membership
test, elimination, projection, and many more. More precisely, the
Gröbner bases framework makes these operations \emph{computationally
accessible} allowing to perform actual computations on ideals.

\subsection{Comparing border bases and Gröbner bases}
\label{sec:comp-bord-bases}

Unfortunately,
Gröbner bases are not always well suited to perform theses operations,
in particular when the actual ideals under consideration are inferred
from measured data. In fact, Gröbner bases do not react smoothly to
small error in the input data.
Border bases are a natural
generalization of Gröbner bases that are believed to 
 deform smoothly in the input (cf. \cite{kreuzer2008dbb}) and that have been used for computations with numerical data (cf., e.g., \cite{HKPP2009, abbott2008sbb}). 

In \cite{HKPP2009} a numerically stable version of the Buchberger\textendash Möller
algorithm has been derived which heavily relies on the usage of border
bases in order to ensure numerical stability. An alternative algorithm
based on border bases with similar stability behavior was investigated
in \cite{abbott2008sbb}. Both algorithms calculate bases of ideals
of points which can be used to derive polynomial models from empirical
datasets. 

It is well known that every Gröbner basis with respect to a degree-compatible
term ordering can be extended to a border basis (see \cite[p.\ 281ff]{kehrein2006cbb})
but not every border basis is an extension of a Gröbner
basis. Moreover, not every order ideal $\oi$ supports a border basis
even if it has the right
cardinality. An example illustrating these
two cases is presented in \cite[Example 6]{kehrein2006cbb}. 

While the border basis algorithm in \cite{kehrein2006cbb}, which is a
specification of Mourrain's generic algorithm \cite{mourrain1999ncn}, allows
for computing border bases of zero-dimensional ideals for order ideals
supported by a degree-compatible term ordering, it falls short to
provide a border basis for more general order ideals: The computed
border basis is supported by a reduced Gröbner basis. The alternative
algorithm presented in \cite[Proposition 5]{kehrein2006cbb} which
can \emph{potentially} compute arbitrary border bases requires
the \emph{a priori} knowledge of the order ideal that might support
a border basis so that the order ideal has to be guessed in advance.
As we cannot expect this prior knowledge, the algorithm does not solve
the problem of characterizing those order ideals for which a border
basis does exist. Further, as pointed out in \cite[p.\ 284]{kehrein2006cbb},
the basis transformation approach of this algorithm is unsatisfactory
as it significantly relies on Gröbner basis computations.

It can be favorable though to be able to compute border bases for
order ideals that do not necessarily stem from a degree-compatible
term ordering.  It has been an open question to characterize those order
ideals of a given zero-dimensional ideal $I$ that support a border basis.
We will answer this question for degree-compatible order ideals by
resorting to polyhedral combinatorics, mentioned in
\autoref{sec:polynomial-method}.  The restriction
to degree-compatible order ideals is due to the design of the algorithm
to proceed degree-wise and the authors believe that this restriction
can be overcome as well.

The preference of border bases over Gröbner bases partly arises from
the iterative generation of linear syzygies, inherent in the border
basis algorithm, which allows for successively approximating the basis degree-by-degree. 
Moreover, the border basis algorithm is a linear algebra algorithm with a tiny exception,
\hyperref[alg:finalReduction]{the final step}, which does not contribute
to the inherent complexity though as
its running time is polynomial in the input size.

\subsection{Polynomial method}
\label{sec:polynomial-method}

In discrete mathematics and combinatorial optimization, polynomial
systems have been used to formulate combinatorial problems such as
the graph coloring problem, the stable set problem, and the matching
problem (see e.g., \cite{deLoera2009ecps} for an extensive list of references).
 This well-known method, which Alon referred to as the \emph{polynomial
method} (cf.\ \cite{alon1999PolyMethod,alon1996PolyMethod}) recently
regained strong interest. This emphasizes the alternative view of the border basis algorithm as a proof system which successively uncovers hidden information by making it explicit.  In \cite[Section 2.3]{deLoera2009cpe} and \cite{deLoera2009ecps,deLoera2008hna}
infeasibility of certain combinatorial problems, e.g., 3-colorability
of graphs is established using Hilbert's (complex) Nullstellensatz
and the authors provide an algorithm \textalgorithm{NulLA} to establish infeasibility
by using a linear relaxation. The core of the algorithm is identical
to the \hyperref[alg:stableSpan]{$L$-stable span procedure} used in the border basis algorithm,
which intimately links both procedures. The difference is of
a technical nature: whereas \textalgorithm{NulLA} establishes
infeasibility,
\hyperref[alg:borderBasis]{the classical border basis algorithm}
as presented in
\cite{kehrein2006cbb} computes the actual border bases of the ideal.
Another recent link between border bases and the Sherali\textendash Adams closure
(introduced in \cite{SA}), a method for convexification, was investigated in \cite{SP20092}.
The authors show that the Sherali\textendash Adams procedure can be understood
as a weaker version of the $L$-stable span procedure, which is an
essential part of the border basis algorithm. As a consequence a new
tighter relaxation than the Sherali\textendash Adams closure is derived by exploiting
the stronger $L$-stable span procedure.

\subsection{Applications of border bases}
\label{sec:appl-bord-bases}

Surprisingly, it turns out that there are deep connections
to other mathematical disciplines, and border bases represent the combinatorial
structure of the ideal under consideration in a canonical way. Although
the use of border basis as a concise framework is quite recent,
introduced in \cite{kehrein2005cbb,kehrein2006cbb,kehrein14asv} and \cite[Section 6.4]{kreuzer2005cca},
the concept of border basis
is rather old and appeared in different fields of mathematics
including computer algebra, discrete optimization, logic, and cryptography
(independently) under different names. Originally, border bases were
introduced in computer algebra as a generalization of Gröbner bases
to address numerical instabilities. Border bases have been successfully
used since for solving zero-dimensional systems of polynomial equations
(see, e.g., \cite{auzinger1988eac,moller1993sae,mourrain1999ncn}),
which in particular include those with solutions in 0/1 and thus a
large variety of combinatorial problems. 

Border bases have been also used to solve sparse quadratic systems
of equations thus giving rise to applications in cryptography in
a natural way. Such systems arise from crypto systems (such as AES,
BES, HFE, DES, CTC variants, etc.) when rewriting the S-boxes as polynomial
equations. The celebrated XL, XSL, MutantXL attacks, which are based
on relinearization methods, are essentially equivalent to the reformulation-linearization-technique
(RLT) of Sherali and Adams \cite{SA} and use a version of the Nullstellensatz
to break ciphers. In fact the XL algorithm (see e.g., \cite{kreuzerAlgAttack,courtois2000eas})
in its classical form is actually identical to a level $d$ Sherali\textendash Adams
closure of the associated system and therefore it is a border bases computation
at its core. Motivated by the success of the aforementioned methods,
border bases have also been used in crypto analysis and coding theory, see
\cite{borgesquintana:ams}.

Another application of border bases is the modeling of dynamic systems
from measured data (see e.g., \cite{HKPP2009, kreuzer2009sbb,
  abbott2008sbb}) where the better numerical stability can be
advantageous.

\subsection*{Our contribution}

A perennial problem in various applications
is that the classical computation of a border basis depends
on a degree-compatible term ordering and
hence finds only the border bases supported
by an order ideal induced by the ordering. These border bases do represent only a rather small fraction of all possible border bases (see Example~\ref{exa:generic}). Further, in this case the border basis contains a reduced
Gröbner basis, and for example the theoretically nice numerical/symmetry properties
are often lost, as the vector space basis of $\pring/I$ (where
$I$ is a zero-dimensional ideal) might not be optimally suited for
numerical computation.  Moreover, for example when solving large systems
of polynomial equations, it is desirable to \emph{guide} the solution
process by having a vector space basis that does actually have a nice
interpretation attached. These bases often cannot be obtained by a
degree-compatible term ordering. 

Techniques from commutative algebra,
and in particular border bases (and alike) have been of great value
to discrete mathematics and combinatorial optimization. This time
we will proceed the other way around. We will apply the combinatorial optimization
toolbox to the combinatorial structure of zero-dimensional ideals
and we solve the aforementioned problems by \emph{freeing the computation
of a border basis from the term ordering}. Our contribution is the
following:

\begin{description}
\item[Polyhedral characterization of all border bases]

We provide a complete, polyhedral characterization of all degree-compatible border bases
of any zero-dimensional ideal $I$:
we associate an \emph{order ideal polytope
$P$} to \(I\) that characterizes all its degree-compatible border bases.
The integral 0/1 points in the order ideal polytope are in one-to-one correspondence with
degree-compatible order ideals supporting a border basis of $I$.
(An order ideal is degree-compatible if it is a degree-wise complement
of \(I\), see Definition~\ref{def:indepSizes}.)
This explicitly establishes the link between the combinatorial structure
of the basis of the factor space and the structure of the ideal. Whether
an order ideal supports a border basis is solely determined
by the combinatorial structure of the order ideal polytope. A related
result for Gröbner bases of a vanishing ideal of generic points was established in \cite{onn1999cutting}.

\item[Computation of border bases not coming from a term ordering]

We outline an algorithm based on the classical border basis algorithm
as defined in \cite{kehrein2006cbb} and show how to compute border
bases for arbitrary degree-compatible order ideals without relying
on a specific degree-compatible term ordering. Recall that not every
order ideal supports a border basis and \emph{guessing} the order
ideal in advance is permissive. It is also not advised to search for
an admissible order ideal with brute force as the combinatorial structure
of the ideal can be so complicated that it is not even possible to
easily guess a (non-term ordering induced) feasible order ideal, not
even to mention to find one that has a preferred structure (see also the next point).

\item[Finding preferred order ideals]
Our approach is also able to compute order ideals
maximizing a prespecified \emph{preference function}
(and compute their border bases).
We will show that computing a preference-optimal order ideal supporting a
border basis is \NP-hard in general, which is surprising, as choosing the
order ideal is merely a basis transformation.
The \NP-hardness does not come from the hardness of
computing the \hyperref[alg:stableSpan]{\(L\)-stabilized span},
as the problem remains \NP-hard in cases
where the \(L\)-stabilized span is small enough to be determined efficiently.

\item[Computational feasibility]

We provide computational tests that demonstrate the feasibility of
our method.

\item[Applications: The structure of ideals and counting of border bases]

Having the order ideal polytope available for a zero-dimensional ideal
$I$, it is possible to examine the structure of the ideal based on
its border bases. One straight-forward application is \emph{counting}
the number of degree-compatible border bases for a zero-dimensional
ideal $I$.
\end{description}

Order ideals and determining those with maximum score do appear in a very
natural way in combinatorial optimization as the so called \emph{maximum weight
  closure problem} (cf.\ e.g., \cite{picard1976maximal}) and they have a
variety of applications, e.g., in open-pit mining where any feasible
production plan is indeed an order ideal; clearly we can only mine the lower
levels \emph{after} having mined the upper ones. A good survey as well as an
introduction to the problem can be found in
\cite{hochbaum2000performance}.

Other, more involved applications might arise,
e.g., in computational biology where the structure of boolean network
is inferred from the Gröbner fan. Due to their higher numerical stability,
border bases might be a better choice here and the order ideal polytope
allows for a detailed study of the underlying combinatorial structure
of the networks. 

\subsection*{Outline}

We start with the necessary preliminaries in \autoref{sec:Preliminaries}
and recall the classical border basis algorithm in \autoref{sub:classicBB}.
In \autoref{sec:The-order-ideal} we introduce the order ideal
polytope and establish the one-to-one correspondence between the 0/1
points of this polytope and degree-compatible border bases. We also
derive an equivalent characterization that is better suited for actual
computations. In \autoref{sec:Computing-border-bases} we then
show how the results from \autoref{sec:The-order-ideal} can be
used to compute admissible degree-compatible border bases using a
preference function for the structure of the order ideal.
We also establish the \NP-hardness of the optimization variant of this
problem.
As an immediate
application, we also outline how to count all admissible border bases
of a zero-dimensional ideal.  We conclude
with computational results in \autoref{sec:Computational-results}
and a few final remarks in \autoref{sec:Concluding-remarks}.

\section{\label{sec:Preliminaries}Preliminaries}

We consider a polynomial ring $\pring$ over the field $K$ with indeterminates
$\mathbb{X}=\{x_{1},\dots,x_{n}\}$. For convenience we define $x^{m}\coloneqq \prod_{j\in[n]}x_{j}^{m_{j}}$
for $m\in\N^{n}$ and let $\T^{n}\coloneqq \{x^{m}\mid m\in\N^{n}\}$
be the \emph{monoid of terms}. For any $d\in\N$ we let $\T_{\leq d}^{n}\coloneqq \{x^{m}\in\T^{n}\mid{\lVert m \rVert}_{1}\leq d\}$
be the set of monomials of total degree at most $d$. For a polynomial
$p\in\pring$ with $p=\sum_{i=1}^{l}a_{i}x^{m_{i}}$ we define the
\emph{support of $p$} to be $\supp(p)\coloneqq \{x^{m_{i}}\mid i\in[l]\}$
and similarly, for a set of polynomials $P\subseteq\pring$ we define
the \emph{support of $P$} to be $\supp(P)\coloneqq \bigcup_{p\in P}\supp(p)$.
Given a term ordering $\sigma$, the \emph{leading term} $\LT_{\sigma}(p)$
of a polynomial $p$ is $\LT_{\sigma}(p)\coloneqq t$ with $t\in\supp(p)$
such that for all $t'\in\supp(P)\setminus\{t\}$ we have $t>_{\sigma}t'$;
the \emph{leading coefficient} $\LC_{\sigma}(p)$ of $p$ is the
coefficient of $\LT_{\sigma}(p)$. We drop the index \(\sigma\) if
the ordering is clear from the context. Recall that the \emph{degree
of a polynomial} $p\in\pring$ is $\deg(p)\coloneqq \max_{x^{m}\in\supp(p)}{\lVert m \rVert}_{1}$.
The \emph{leading form} $\LF(p)$ of a polynomial $p=\sum_{i=1}^{l}a_{i}x^{m_{i}}\in\pring$
is defined to be $\LF(p)=\sum_{i=1,{\lVert m_{i} \rVert}=d}^{l}a_{i}x^{m_{i}}$
where $d=\deg(p)$, i.e., we single out the part with maximum degree.
Both $\LF$ and $\LT$ generalize to sets in the obvious way, i.e.,
for a set of polynomials $P$ we define $\LF(P)\coloneqq \{\LF(p)\mid p\in P\}$
and $\LT(P)\coloneqq \{\LT(p)\mid p\in P\}$.

In the following we
will frequently switch between considering polynomials $M=\{p_{1},\dots,p_{s}\}$,
the generated ideal, and the generated vector space whose coordinates
are indexed by the monomials in the support of $M$. We denote the ideal
generated by $M$ as $\igen{M}$ and the vector space generated by
$M$ as $\vgen{M}$. For $n\in\N$ we define $[n] \coloneqq \{1,\dots,n\}$.
All other notation is standard as to be found in \cite{cox2007iva,kreuzer2000cca};
we have chosen the border basis specific notation to be similar to the one
in \cite{kehrein2006cbb} where the border basis algorithm was first
introduced in its current form.

Central to our discussion will be the notion of an order ideal:
\begin{defn}
Let $\oi$ be a finite subset of $\T^{n}$. If for all $t\in\oi$
and $t'\in\T^{n}$ such that $t'\mid t$ we have $t'\in\oi$, i.e.,
$\oi$ is closed under factors, then we call $\oi$ an \emph{order
ideal}. Furthermore, the \emph{border $\partial\oi$} of a non-empty
order ideal $\oi$ is the set of terms $\partial\oi\coloneqq \{x_{j}t\mid j\in[n],t\in\oi\}\setminus\oi$.
As an exception,
we set $\partial\emptyset\coloneqq \{1\}$ for the empty order ideal.
\end{defn}
Recall that an ideal $I\subseteq\pring$ is zero-dimensional, if and only if $\pring/I$
is finite dimensional. The $\oi$-border basis of a zero-dimensional
ideal $I$ is a special set of polynomials:
\begin{defn}\label{def:bbasis}
Let $\oi=\{t_{1},\dots,t_{\mu}\}$ be an order ideal with border $\partial\oi=\{b_{1},\dots,b_{\nu}\}$.
Further let $\mathcal{G}=\{g_{1},\dots,g_{\nu}\}\subseteq\pring$
be a (finite) set of polynomials and let $I\subseteq\pring$ be a
zero-dimensional ideal. Then the set $\mathcal{G}$ is an \emph{$\oi$-border
basis} of \(I\) if:
\begin{enumerate}
\item \label{enu:bbform} the polynomials in $\mathcal{G}$ have the form $g_{j}=b_{j}-\sum_{i=1}^{\mu}\alpha_{ij}t_{i}$
for $j\in[\nu]$ and $\alpha_{ij}\in K$;
\item $\igen{\mathcal{G}}=I$;
\item $\pring=I\oplus\vgen{\oi}$ as vector spaces.
\end{enumerate}
If there exists an $\oi$-border basis of $I$ then the order
ideal $\oi$ \emph{supports a border basis} \emph{of $I$}.
\end{defn}
Note that the condition $\igen{\mathcal{G}}=I$ is
a consequence of $\mathcal{G}\subseteq I$, the particular form of the
elements in $\mathcal{G}$, and $\pring=I\oplus\vgen{\oi}$.
More precisely, \(\igen{\mathcal{G}} + \vgen{\oi}\) is closed under
multiplication by the \(x_{i}\), and hence it is an ideal.
As it contains \(1\), it must be the whole ring, and hence
$\igen{\mathcal{G}}=I$ by the modular law.
See \cite[Proposition 4.4.2]{kehrein14asv} for another proof.

In particular, an order ideal
$\oi$ supports an $\oi$-border basis of $I$
if and only if
$\pring = I \oplus \vgen{\oi}$.
Moreover, for any given order ideal
$\oi$ and ideal $I$ the $\oi$-border basis of $I$ is unique as
$b_{j}$ has a unique representation in $\pring=I\oplus\vgen{\oi}$
for all $j\in[\nu]$. Furthermore, as $\pring=I\oplus\vgen{\oi}$
it follows that $\lvert \oi \rvert = \dim \vgen{\oi} $ is invariant for all
choices of $\oi$.
The requirement for $I$ being zero-dimensional is a consequence of
this condition as well, as $\partial\oi$ (and in consequence $\oi$)
as part of the output of our computation should be finite. Clearly,
as a vector space, $I$ has a degree filtration, i.e., $I=\bigcup_{i\in\N}I^{\leq i}$
where $I^{\leq i}\coloneqq \{ p\in I\mid\deg(p)\leq i \}$. For a set of
monomials $\oi$ we define $\oi^{=i}\coloneqq \{m\in\oi\mid\deg(m)=i\}$
(and similarly for $\leq$ instead of \(=\)). In the following we consider special
types of order ideals, i.e., those that essentially preserve the filtration: 
\begin{defn}
\label{def:indepSizes}Let $I\subseteq\pring$ be a zero-dimensional
ideal and let $\oi\subseteq\T^{n}$ be an order ideal. We say that
$\oi$ is \emph{degree-compatible} (to $I$) if
\[
\dim \vgen{\oi^{=i}} =\dim \vgen{\T_{=i}^{n}} -\dim \frac{I^{\leq i}}{I^{\leq i-1}} \]
for all $i\in\N$. 
\end{defn}
Thus, the $\oi$-border basis of a zero-dimensional ideal $I$ with
respect to any degree-compatible order ideal $\oi$ has a pre-determined
size for each degree $i\in\N$. Intuitively, the degree-compatible
order ideals are precisely those that correspond to degree-compatible
orderings on the monomials. The important difference is that the orderings
do not have to be term orderings. The definition above only requires
\emph{local} compatibility with multiplication as $\oi$ is an order
ideal and thus downwardly closed and degree-compatible, i.e., if $p,q$
are polynomials and $\deg(p)<\deg(q)$ then $p\leq q$. The following
example shows that the requirements of Definition \ref{def:indepSizes}
are not automatically satisfied by all order ideals $\oi$. 
We also give an example of a degree-compatible ideal
not coming from a term ordering.

Recall that the order ideal associated to a Gröbner basis of an ideal
consists of all monomials not divisible by
any leading term in the Gröbner basis.
\begin{example}[A non-degree-compatible order ideal]
With the degree lexicographic ordering on two variables,
the set of polynomials
\[ x_1^3 + x_1 x_2, x_1^2 x_2, x_1 x_2^2, x_2^3\]
is a Gröbner basis of an ideal with associated order ideal
consisting of all monomials of degree at most \(2\):
\[ \{ 1, x_1, x_2, x_1^2, x_1 x_2, x_2^2 \}. \]
The element \(x_1^3 + x_1 x_2\) of the ideal enables us
to replace \(x_1 x_2\) with \(x_1^3\) and thus obtain another
order ideal of the ideal, which is not degree-compatible:
\[ \{ 1, x_1, x_2, x_1^2, x_1^3, x_2^2 \}. \]
\end{example}

\begin{example}[A degree-compatible order ideal not coming from a term ordering]
The homogeneous ideal having a Gröbner basis
in the degree lexicographic term ordering
\[ x_2^2 + x_1 x_2 + x_1^2, x_1x_2^2,x_2^4 \]
has the associated order ideal
\[ \{ 1, x_1, x_2, x_1 x_2, x_2^2, x_2^3 \}. \]
Via a basis change, we obtain another order ideal
\[  \{ 1, x_1, x_2, x_1^2, x_2^2, x_2^3 \}, \]
which cannot come from a term ordering.
The reason is that every Gröbner basis of the ideal must contain
a polynomial with degree-two leading term, which therefore
must be a multiple of the first generator
\(x_2^2 + x_1 x_2 + x_1^2\)
modulo higher-degree terms.
The degree-two leading term cannot be \(x_1 x_2\) because
if e.g., \(x_1 < x_2\) then the term \(x_2^2\) is larger then \(x_1 x_2\).
So the leading term must be either \(x_1^2\) or \(x_2^2\),
which excludes the order ideal.
\end{example}

\begin{example}[Generic ideal]\label{exa:generic}
  Let \(k\) and \(n\) be positive integers and
  let \({\{a_{ij}\}}_{i \in [n], j \in [k]}\) be
  algebraically independent real numbers over \(\Q\).
  Let \(I\) be the ideal of polynomials in the variables \(x_1,\dotsc,x_n\)
  which are zero on the points \((a_{1j}, \dotsc, a_{nj})\)
  for \(j \in [k]\).
  Thus, the ideal is zero-dimensional, and \(\pring / I\) has dimension \(k\).

  Every \(k\) distinct monomials form a complementary basis of \(I\),
  since they are linearly independent on the \(k\) points
  \((a_{1j}, \dotsc, a_{nj})\).
  An equivalent formulation of linear independence is that
  the determinant of the matrix formed by the values of
  the monomials on these points is non-zero.
  The determinant is indeed non-zero,
  as it is a non-trivial polynomial of the algebraic independent \(a_{ij}\)
  with integer coefficients.
  
  In particular, every order ideal of size $k$ is an order ideal of \(I\).
  The degree-compatible order ideals are the ones
  where the monomials have the least possible degree,
  i.e., consisting of
  all monomials of degree less than \(l\) and
  \(k - \binom{n+l-1}{l-1}\) monomials of degree \(l\),
  where \(l\) is the smallest non-negative integer satisfying
  \(k \leq \binom{n+l}{l}\), i.e.,
  there are at least \(k\) monomials of degree at most \(l\).
\end{example}

\subsection{\label{sub:classicBB}Computing border bases for term ordering induced
$\oi$}

Without proofs,
we recall the classical border basis algorithm
introduced in \cite{kehrein2006cbb} as it will serve as a basis for
our algorithm.  The interested reader is referred to
\cite{kehrein14asv,kehrein2005cbb} for a general introduction to
border bases and to \cite{kehrein2006cbb} in particular for an introduction
of the border basis algorithm.

The classical border basis algorithm
calculates border bases of zero-dimensional ideals with
respect to an order ideal $\oi$ which is induced by a degree-compatible
term ordering $\sigma$.

First, the border basis algorithm heavily relies
on the following \emph{neighborhood extension}:
\begin{defn}
\label{def:lstabSpan}
(cf.\ \cite[Definition 7.1 and paragraph preceding Proposition 13]{kehrein2006cbb})
Let $V$ be a vector space. We define the \emph{neighborhood
extension of $V$} to be
\[V^{+}\coloneqq V+Vx_{1}+\dots+Vx_{n}.\]

For a finite set \(W\) of polynomials,
its \emph{neighborhood extension} is
\[ W^{+}=W\cup Wx_{1}\cup\dots\cup Wx_{n}. \]
\end{defn}
Note that for a given set of polynomials $W$ such that $\vgen{W}=V$
we have $\vgen{W^{+}}=\vgen{W}^{+}=V^{+}$
as multiplication with $x_{i}$ is a $K$-vector space homomorphism.
It thus suffices to perform the neighborhood extension on a set of
generators $W$ of $V$. 

Let $F$ be a finite set of polynomials and let $L\subseteq\T^{n}$
be an order ideal, then $F\cap \vgen{L} = \{f\in F\mid\supp(f)\subseteq L\}$,
i.e., $F\cap \vgen{L}$ contains only those polynomials that lie in the vector
space generated by $L$. Clearly, for $L = \T^n_{\leq d}$ we have $\vgen{F} \cap \vgen{L} = \vgen{F}^{\leq d}$. Using the neighborhood extension we can define:
\begin{defn}(cf.\ \cite[Definition 10]{kehrein2006cbb})
Let $L$ be an order ideal and let $F$ be a finite
set of polynomials such that $\supp(F)\subseteq L$.
The set \(F\) is \emph{$L$-stabilized} if
$\vgen{F^{+}} \cap \vgen{L} = \vgen{F}$.
The \emph{$L$-stable span} \(F_L\) of \(F\) is
the smallest vector space \(G\) containing \(F\) satisfying \(G^+ \cap \vgen{L} = G\).
\end{defn}
A straightforward construction of the \(L\)-stable span of \(F\) is
to inductively define the following increasing sequence of vector spaces:
\[
F_{0} \coloneqq \vgen{F} \quad\text{and}\quad
F_{k+1} \coloneqq F_{k}^{+} \cap \vgen{L}\text{ for \(k>0\)}.
\]
The union $\bigcup_{k\geq0} F_{k}$ is the $L$-stable span \(F_{L}\) of \(F\).

The set $L$ represents our \emph{computational universe} and we will
be in particular concerned with finite sets $L\subseteq\T^{n}$. Note
that $L$-stability is a property of the vector space and does not
depend on the basis:
\begin{rem}
\label{rem:stableForVR}
The \(L\)-stable span of a finite set \(F\) depends only on
the generated vector space \(\vgen{F}\), as \(\vgen{F^{+}} = \vgen{F}^{+}\). 
\end{rem}

In the following we will explain how the $L$-stable span can be computed
explicitly for $L=\T_{\leq d}^{n}$. We will use a modified version
of Gaussian elimination as a tool, which allows us to extend a given
basis $V$ with a set $W$ as described in the following:
\begin{lem}
\label{lem:borderGauss}\cite[Lemma 12]{kehrein2006cbb} Let $V=\{v_{1},\dots,v_{r}\}\subseteq \pring \setminus\{0\}$
be a finite set of polynomials such that $\LT(v_{i})\neq\LT(v_{j})$
whenever $i,j\in[r]$ with $i\neq j$ and $\LC(v_{i})=1$ for all
$i\in[r]$. Further let $G=\{g_{1},\dots,g_{s}\}$ be a finite set
of polynomials. Then Algorithm \ref{alg:gaussEl} computes a finite
set of polynomials $W\subseteq \pring$ with $\LC(w)=1$ for all $w\in W$,
$\LT(u_{1})\neq\LT(u_{2})$ for any distinct $u_{1},u_{2}\in V\cup W$,
and $\vgen{V\cup W}=\vgen{V\cup G}$. ($V$, $W$
may be empty.)\end{lem}
\begin{algorithm}{Gaussian Elimination for polynomials}{\gausEl}%
{$V$, $G$ as in Lemma \ref{lem:borderGauss}.}%
{$W\subseteq \pring$ as in Lemma \ref{lem:borderGauss}.}
\label{alg:gaussEl}

\item Let $H\coloneqq G$ and $\eta\coloneqq 0$.
\item \label{enu:step2}If $H=\emptyset$ then return $W\coloneqq \{v_{r+1},\dots,v_{r+\eta}\}$
and stop.
\item Choose $f\in H$ and remove it from $H$. Let $i\coloneqq 1$.
\item \label{enu:step4}If $f=0$ or $i>r+\eta$ then go to step (\ref{enu:step7}).
\item If $\LT(f)=\LT(v_{i})$ then replace $f$ with $f - \LC(f) \cdot v_{i}$.
Set $i\coloneqq 1$ and go to step (\ref{enu:step4}).
\item Set $i\coloneqq i+1$. Go to step (\ref{enu:step4}).
\item \label{enu:step7}If $f\neq0$ then put $\eta\coloneqq \eta+1$ and let $v_{r+\eta}\coloneqq f/\LC(f)$.
Go to step (\ref{enu:step2}).
\end{algorithm}

We can now compute the $L$-stable span using the Gaussian elimination
algorithm \ref{alg:gaussEl}:
\begin{lem}
\cite[Proposition 13]{kehrein2006cbb}\label{lem:stableSpan} Let
$L=\T_{\leq d}^n$ and
$F \subseteq \pring$ be a finite set of polynomials supported on \(L\).
Then Algorithm \ref{alg:stableSpan} computes a vector space basis
$V$ of $F_{L}$ with pairwise different leading terms.\end{lem}
\begin{algorithm}{$L$-stable span computation}{\lstabspan}%
{$F$, $L$ as in Lemma \ref{lem:stableSpan}.}%
{$V$ as in Lemma \ref{lem:stableSpan}.}
\label{alg:stableSpan}
\item $V\coloneqq \gausEl(\emptyset,F)$.
\item \label{enu:step2StabSpan}$W'\coloneqq \gausEl(V,V^{+} \setminus V)$.
\item\label{enu:step3StabSpan}
  $W\coloneqq \{w\in W'\mid\supp(w) \subseteq L\} = \{ w \in W' \mid
  \deg(w) \leq d \}$.
\item If $\lvert W \rvert >0$ set $V \coloneqq V \cup W$ and go to step
  (\ref{enu:step2StabSpan}).
\item Return \(V\).
\end{algorithm}
The last ingredient that we need in order to formulate the border
basis algorithm is the final reduction algorithm. This algorithm basically
transforms the output of the border basis algorithm to bring it into
the desired form of a border basis by applying linear algebra steps
only. It interreduces the elements so that they only have support
in the leading term and $\oi$. 
\begin{lem}
\cite[Proposition 17]{kehrein2006cbb}\label{lem:finalRed} Let $F$
be a system of generators of a zero-dimensional ideal $I$, let $L$
be an order ideal, and let $V$ be a vector space basis of $F_{L}$ with pairwise
different leading terms and $\oi\coloneqq L\setminus\LT(V)$ such
that $\partial\oi\subseteq L$.
Then Algorithm \ref{alg:finalReduction} computes the $\oi$-border
basis $\mathcal{G}$ of $I$. \end{lem}
It is easy to see that $\vgen{L}=F_{L}\oplus\vgen{\oi}$,
so the algorithm will do a form of Gaussian elimination
on the basis \(V\) to obtain a new basis \(V_{R}\)
with all non-leading terms supported on \(\oi\).
Of course, this new basis will contain an \(\oi\)-border basis.

\begin{algorithm}{Final Reduction Algorithm}{\finalRed}%
{$V$, $\oi$ as in Lemma \ref{lem:finalRed}.}%
{$\mathcal{G}$ as in Lemma \ref{lem:finalRed}.}
\label{alg:finalReduction}

\item Let $V_{R}\coloneqq \emptyset$.
\item \label{enu:finRedStep2}If $V=\emptyset$ then go to step (\ref{enu:finRedStop}).
\item Let $v\in V$ such that $v$ has minimal leading term. Put $V\coloneqq V\setminus\{v\}$.
\item Let $H\coloneqq \supp(v)\setminus(\LT(v)\cup\oi)$.
\item If $H=\emptyset$ then append $v/\LC(v)$ to $V_{R}$ and go to step
(\ref{enu:finRedStep2}).
\item For each $h\in H$ choose $w_{h}\in V_{R}$ and $c_{h}\in K$ such
that $\LT(w_{h})=h$ and $h\notin\supp(v-c_{h}w_{h})$.
\item Set $v\coloneqq v-\sum_{h\in H}c_{h}w_{h}$, append $v/\LC(v)$ to $V_{R}$,
and go to step (\ref{enu:finRedStep2}).
\item\label{enu:finRedStop}
  Return $\mathcal{G}\coloneqq \{ v \in V_{R} \mid \LT(v) \in \partial\oi \}$.
\end{algorithm}

We will now formulate the border basis algorithm.
\begin{prop}
\cite[Proposition 18]{kehrein2006cbb} \label{pro:bbasis}Let $F \subseteq \pring$
be a finite set of polynomials that generates a zero-dimensional ideal
$I=\igen{F}$. Then Algorithm \ref{alg:borderBasis} computes
the $\oi$-border basis $\mathcal{G}$ of $I$.\end{prop}
\begin{algorithm}{Border basis algorithm}{\bbasis}%
{$F$ as in Proposition \ref{pro:bbasis}.}%
{$\mathcal{G}$ as in Proposition \ref{pro:bbasis}.}
\label{alg:borderBasis}

\item Let $d\coloneqq \max_{f\in F} \deg(f) $.
\item \label{enu:bbasisStep2}$V=\{v_{1},\dots,v_{r}\}\coloneqq \lstabspan(F,\T_{\leq d}^{n})$.
\item Let $\oi\coloneqq \T_{\leq d}^{n}\setminus\{\LT(v_{1}),\dots,\LT(v_{r})\}$.
\item \label{enu:bbasisStop}If $\partial\oi\nsubseteq \T_{\leq d}^{n}$ then set $d\coloneqq d+1$
 and go to step (\ref{enu:bbasisStep2}).
\item Return $\mathcal{G}\coloneqq \finalRed(V,\oi)$.
\end{algorithm}
It is worthwhile to note that step~(\ref{enu:bbasisStop})
in Algorithm \ref{alg:borderBasis}
is essentially for testing if the $L$-stable span is large enough
to support an $\oi$-border basis.

The rationale for searching such a large span is due to
the following proposition
that serves as a stopping criterion.
It is obvious from its Corollary~\ref{cor:idealApprox}
that the span is indeed minimal in the sense that
it is the smallest span that contains a degree-compatible order ideal
that supports a border basis (and thus all such degree-compatible order ideals).

\begin{prop}
\label{pro:bbApproximation}
\cite[Proposition 16]{kehrein2006cbb}
Let $L$ be an order ideal.
Further let $\tilde{I}$ be an \(L\)-stabilized generating vector subspace of a zero-dimensional
ideal $I\subseteq\pring$, i.e., ${\tilde{I}}^{+}\cap \vgen{L}=\tilde{I}$
and $\igen{\tilde{I}}=I$. If $\oi$ is an
order ideal such that $\vgen{L}=\tilde{I}\oplus\vgen{\oi}$ and $\partial\oi\subseteq L$
then $\oi$ supports a border basis of $I$.
\end{prop}

We obtain the following corollary:

\begin{corollary}
  \label{cor:idealApprox}
  Let \(\tilde{I}\) be an \(\T^{n}_{\leq d}\)-stabilized vector space
  satisfying \(\tilde{I} + \vgen{\T^{n}_{\leq d-1}} = \vgen{\T^{n}_{\leq d}}\).
  Then \(\igen{\tilde{I}} \cap \vgen{\T^{n}_{\leq d}} = \tilde{I}\).
\begin{proof}
We apply Proposition~\ref{pro:bbApproximation} with the choice
\(L \coloneqq \T^{n}_{\leq d}\),
\(I \coloneqq \igen{\tilde{I}}\)
and
\(\oi \coloneqq \T^{n}_{\leq d} \setminus \LT(\tilde{I})\)
where the leading terms are
with respect to any degree-compatible term ordering.
Clearly,
\(\vgen{\T^{n}_{\leq d}}= \tilde{I} \oplus \vgen{\oi}\).
The condition
\(\tilde{I} + \vgen{\T^{n}_{\leq d-1}} = \vgen{\T^{n}_{\leq d}}\)
ensures that
\(\oi\) consists of monomials of degree less than \(d\),
so \(\partial\oi \subseteq \T^{n}_{\leq d}\).
Hence the proposition applies,
and we obtain \(\pring = I \oplus \vgen{\oi}\).
Together with
\(\vgen{\T^{n}_{\leq d}}= \tilde{I} \oplus \vgen{\oi}\)
this gives
\(I \cap \vgen{\T^{n}_{\leq d}} = \tilde{I}\).
\end{proof}
\end{corollary}

The border basis algorithm decomposes
into two components. The first one is the calculation of the $L$-stable
span for $L=\T_{\leq d}^{n}$ for sufficiently large $d$ (this is the main `work')
and the second component is the extraction of a border basis
via the \hyperref[alg:finalReduction]{final reduction algorithm}.
Effectively,
after using \emph{any} degree-compatible term ordering in order to
compute the $L$-stable span one can choose a different ordering with
respect to which the basis can be transformed. This approach is very
much in spirit of the FGLM algorithm, Steinitz's exchange lemma, and
related basis transformation procedures. In Remark \ref{rem:stableForVR}
it is shown that $L$-stability does not depend on the basis and thus
any sensible basis transformation that results in an admissible order
ideal $\oi$ is allowed. In the following we will characterize all
admissible order ideals.

\section{\label{sec:The-order-ideal}The order ideal polytope}

We will now introduce the \emph{order ideal polytope} $P$ (a 0/1 polytope)  
that characterizes all order ideals that
support a border basis (for a given zero-dimensional ideal $I$) in an abstract fashion independent of particular vector space bases for the stable span approximation. Its role will be crucial for the \hyperref[alg:borderBasisNew]{later computation
  (Algorithm~\ref*{alg:borderBasisNew})} of
border bases for general degree-compatible
order ideals.
In the first subsection,
we introduce the polytope in an abstract, invariant way
highlighting its main property that
its integral points are in bijection with degree-compatible order ideals
supporting a border basis
(Theorem~\ref{thm:polytope=order-ideal}).
In the second subsection,
we give a more direct reformulation targeted to actual computations.

Given $d\in\N$ in advance for which Algorithm \ref{alg:borderBasis}
stops, $\lstabspan$ computes the actual border basis (and some unused
extra polynomials). The computation of the border basis is performed
by Gaussian elimination applied to the matrices obtained within $\lstabspan$.
The order ideal for which we effectively compute the border basis
is solely determined by the pivoting rule when doing the elimination
step. Degree-compatible term orderings ensure that we obtain an order
ideal for which a border basis exists (given the correct $d$). If
we would now pivot arbitrarily, which effectively means permuting
columns, it is not clear that, first, the resulting ideal is an order
ideal and, second, that it does actually support a border basis. If
we remain in the setting of degree-compatible order ideals (which
means that we only permute columns of monomials with same degree)
then we face the combinatorial problem that when permuting a column
such that it will end up being an element of $\oi$ we also have to
ensure that all its divisors will also end up in $\oi$. On the other
hand, due to the degree-compatibility constraint (see Definition \ref{def:indepSizes})
we must choose an exactly determined amount of elements for each degree
that will end up in $\oi$. We will show now how to transform this
combinatorial problem of \emph{faithful pivoting} into a polyhedral
setting. We obtain a $0/1$ polytope $P$, the \emph{order ideal polytope}
that characterizes all admissible degree-compatible order ideals that
support a border basis for the ideal at hand.

\subsection{Theoretical point of view}
\label{sec:theor-viewp}

\begin{defn}
Let \(I\) be a zero-dimensional ideal.
Its \emph{order ideal polytope} $P(I)$
is given by the system of inequalities
in \autoref{fig:OrderIdealPolytope2}.

\begin{figure}
\begin{empheq}[box=\fbox]{align}
  z_{m_{1}} &\geq z_{m_{2}} \quad
  \forall m_{1},m_{2}\in \T^n\colon m_{1}\mid m_{2}\label{eq:orderIdeal2}\\
  \sum_{m\in \T^n _{=i}}z_{m} &= \dim\vgen{\T^n_{=i}} -
  \dim \left. I^{\leq i} \middle/ I^{\leq i-1} \right. \quad
  \forall i \label{eq:degreeCompatibility2}\\
  \label{eq:steinitz2}
  \sum_{m\in U}z_{m} &\leq
  \begin{multlined}[t]
    \dim\vgen{U \cup \left. I^{\leq i} \middle/ I^{\leq i-1} \right.} -
    \dim \left. I^{\leq i} \middle/ I^{\leq i-1} \right. \\
    \forall i,U \subseteq \T^n_{=i}\colon \lvert U \rvert =
    \dim\vgen{\T^n_{=i}}-\dim \left. I^{\leq i} \middle/ I^{\leq i-1} \right.
  \end{multlined}
  \\
  z_{m} &\in [0,1] \quad \forall m\in \T^n \label{eq:0-1}
\end{empheq}
\caption{\label{fig:OrderIdealPolytope2}Order ideal polytope $P(I)$}
\end{figure}
\end{defn}
The order ideal polytope is actually a finite dimensional polytope
as all the \(z_m\) are \(0\) when the degree of \(m\) is large enough.
Indeed,
for large \(i\), we have
\(\vgen{\T^{=i}} \cong \left. I^{\leq i} \middle/ I^{\leq i-1} \right.\)
and hence \autoref{eq:degreeCompatibility} gives \(z_m = 0\)
for every \(m\) of degree \(i\).

We are ready to relate the order ideal polytope with
order ideals.
From now on,
let \(\Lambda(I)\) denote the set of degree-compatible order ideals
of a zero-dimensional ideal \(I\).
\begin{thm}
  \label{thm:polytope=order-ideal}
  Let \(I\) be a zero-dimensional ideal.
  There is a bijection between
  the set $\Lambda(I)$ of its degree-compatible order ideals and
  the set of integral points of the order ideal polytope of \(I\).
  The bijection is given by
  \[
  \xi \colon z \in P(I) \cap \Z^{\T^n} \mapsto
  \oi(z) \coloneqq \{ m \in \T^n \mid z_{m}= 1\}.
  \]
\begin{proof}
In fact, we will see that the order ideal polytope is defined
exactly to this end.

First,
the integral solutions \(z\) of \autoref{eq:0-1}
are exactly the 0/1 points, i.e.,
the characteristic vectors
of sets of terms $\oi(z)\coloneqq \{ m \in \T^n \mid z_{m}=1 \}$.

Second, it is easy to see that
\autoref{eq:orderIdeal} means that
$\oi(z)$ is indeed an order ideal, as whenever
$m_{1} \mid m_{2}$ and $m_{2}\in\oi(z)$, i.e., $z_{m_{2}}=1$, then
it follows that $z_{m_{1}}=1$, i.e., $m_{1}\in\oi(z)$ as well.

In the third step we will provide an algebraic characterization
of \equationautorefname s \ref{eq:degreeCompatibility2} and \ref{eq:steinitz2}.
Clearly, \autoref{eq:degreeCompatibility2}
can be rewritten to
\[
\lvert {\oi(z)}^{=i} \rvert = \dim\vgen{\T^n_{=i}} -
\dim \left. I^{\leq i} \middle/ I^{\leq i-1} \right. .
\]

We will now show that \autoref{eq:steinitz2} is equivalent
to
\begin{equation}
\label{eq:disjoint}
\vgen{{\oi(z)}^{=i}}\cap\left. I^{\leq i} \middle/ I^{\leq i-1} \right.=\{0\},
\end{equation}
i.e., the image of \({\oi(z)}^{=i}\) is linearly independent in the factor
\(\left. \vgen{\T^n_{=i}} \middle/
  \left( I^{\leq i} \middle/ I^{\leq i-1} \right) \right.\).

Similarly as above,
\autoref{eq:steinitz2}
can be rewritten to
\begin{equation*}
    \left\lvert U \cap {\oi(z)}^{=i} \right\rvert \leq
      \dim\vgen{U \cup \left. I^{\leq i} \middle/ I^{\leq i-1} \right.} -
      \dim \left. I^{\leq i} \middle/ I^{\leq i-1} \right.,
\end{equation*}
i.e., the size of \(U \cap {\oi(z)}^{=i}\) is at most
the dimension of the vector space generated by the image of \(U\) in
the factor
\(\left. \vgen{\T^n_{=i}} \middle/
  \left( I^{\leq i} \middle/ I^{\leq i-1} \right) \right.\).
This is obviously \emph{necessary} for
the image of \({\oi(z)}^{=i}\) to be linearly independent in the factor.
(Here the size of \(U\) does not matter.)

For \emph{sufficiency}
choose $U\coloneqq {\oi(z)}^{=i}$. Then
the dimension of the vector space generated by the image of \({\oi(z)}^{=i}\)
is at least
\(\lvert{\oi(z)}^{=i}\rvert\),
so the image of \({\oi(z)}^{=i}\)  is independent. 

So far we have proved that the integral points of the order ideal polytope
correspond bijectively to order ideals \(\oi\) with the properties
\[
\lvert {\oi(z)}^{=i} \rvert = \dim \vgen{\T^n_{=i}} -
\dim \left. I^{\leq i} \middle/ I^{\leq i-1} \right.
\]
and
\[
\vgen{{\oi(z)}^{=i}} \cap \left. I^{\leq i} \middle/ I^{\leq i-1} \right. = \{0\}
\]
for all \(i\).

Lastly, we will show now that this is equivalent to
$\lvert {\oi(z)}^{=i} \rvert = \dim \vgen{\T^n_{=i}} -
\dim \left. I^{\leq i} \middle/ I^{\leq i-1} \right.$
and $I \oplus \vgen{\oi(z)} = \vgen{\T^n}$ and thus the assertion follows.

By dimensionality it follows, 
$\lvert {\oi(z)}^{=i} \rvert = \dim \vgen{\T^n_{=i}} -
\dim \left. I^{\leq i} \middle/ I^{\leq i-1} \right.$
and
$\vgen{{\oi(z)}^{=i}} \cap \left. I^{\leq i} \middle/ I^{\leq i-1} \right. =
\{0\}$
for all $i$
together are equivalent to
$\left. I^{\leq i} \middle/ I^{\leq i-1} \right.\oplus \vgen{{\oi(z)}^{=i}} =
\vgen{\T^n_{=i}}$
for all $i$.
For brevity,
we will omit the phrase \lq for all \(i\)\rq.
Using the filtration argument,
the latter is equivalent to
$I^{\leq i} \oplus \vgen{{\oi(z)}^{\leq i}} = \vgen{\T^n_{\leq i}}$.
Via a dimension argument on embeddings,
this is further equivalent to
$I \oplus \vgen{\oi(z)} = \vgen{\T^n}$
and
$\dim I^{\leq i} + \dim \vgen{{\oi(z)}^{\leq i}} = \dim\vgen{\T^n_{\leq i}}$.
Finally, filtrating the dimension by degree shows that the latter
is equivalent to
$I \oplus \vgen{\oi(z)} = \vgen{\T^n}$
and
$\dim \left. I^{\leq i} \middle/ I^{\leq i-1} \right. +
\dim \vgen{{\oi(z)}^{=i}} = \dim\vgen{\T^n_{=i}}$.
\end{proof}
\end{thm}

\subsection{Computational point of view}
\label{sec:polytope-comp-viewp}

Throughout this subsection we assume that
$M=\bigcup_{i\in\N}M_{i}\subseteq\pring$
is a finite set of polynomials with degree-filtration $\{M_{i}\mid i\in\N\}$
that generates a zero-dimensional ideal $\igen{M}\subseteq\pring$
such that all $p\in M_{i}$ have degree \(i\).  Furthermore, for
each $i\in\N$ we have an enumeration $M_{i}=\{p_{ij}\mid j\in[k_{i}]\}$
with $k_{i}\in\N$. As seen in \autoref{sub:classicBB}, an important
component of the border basis algorithm is the computation of $L$-stable
spans with respect to some \emph{computational universe} $L\subseteq\T^{n}$.
Computing a border basis with respect to a different order ideal $\oi$
is merely a basis transformation of the vector space obtained from
the $\lstabspan$ procedure.
In the following we assume
that $L$ is of the form $L=\T_{\leq d}^{n}$ where $d\in\N$ is such that Algorithm \ref{alg:borderBasis}
stops.
In view of Remark \ref{rem:stableForVR}
and Definition \ref{def:indepSizes} we have that for a given computational
universe $L$ either all border bases (supported by degree-compatible
order ideals) are contained in $L$ or none.
Furthermore we assume that $M$ is $L$-stabilized and has in particular
a convenient form. Effectively one might want to think of $M$ being
the output of the $\lstabspan$ procedure, which is then brought into
the following reduced form:
\begin{defn}
\label{def:L-canonical-form}
Let $M$ be a finite set of polynomials of degree at most \(\ell\)
for some $\ell\in\N$.
Then \(M\) is in \emph{canonical form}
if the leading term of any element of \(M\) does not occur
in the other elements.
Here we can freely choose leading terms of the polynomials
with the only constraint that they have to be maximal-degree terms.
\end{defn}
We give a visual interpretation of the definition.
The coefficient matrix
 $A\in K^{M\times\T_{\leq \ell}^{n}}$ of \(M\)
is the matrix where
the rows are the elements of \(M\),
the columns are all the monomials of degree at most \(\ell\),
and the entries are the coefficients of the terms in the elements of \(M\).
We use the convention
that for terms $t_{1}$, $t_{2}$ with $\deg(t_{1})>\deg(t_{2})$
we put column $t_{1}$ to the left of column $t_{2}$.
Similarly, we put leading terms of a polynomial to
the left of the other terms.
Now $M$ is in \emph{canonical form} if
the matrix $A$ has
the structure as depicted in \autoref{fig:l-can-form}, i.e., it
consists of degree blocks and each degree block is maximally interreduced. 

\begin{figure}
\[
A=\left(\begin{array}{ccc|c||ccc|c||c||ccc|c}
1 &  & 0 & \star &  &  &  & \star &  &  &  &  & \star\\
 & \ddots &  & \vdots &  & 0 &  & \vdots &  &  & 0 &  & \vdots\\
0 &  & 1 & \star &  &  &  & \star &  &  &  &  & \star\\
\hline  &  &  & 0 & 1 &  & 0 & \star &  &  &  &  & \star\\
 & 0 &  & 0 &  & \ddots &  & \vdots &  &  & 0 &  & \vdots\\
 &  &  & 0 & 0 &  & 1 & \star &  &  &  &  & \star\\
\hline  &  &  &  &  &  &  &  & \ddots & & & & \\
\hline  &  &  & 0 &  &  &  & 0 &  & 1 &  & 0 & \star\\
 & 0 &  & 0 &  & 0 &  & 0 &  &  & \ddots &  & \vdots\\
 &  &  & 0 &  &  &  & 0 &  & 0 &  & 1 & \star\end{array}\right)\]

\caption{\label{fig:l-can-form}canonical form}

\end{figure}

The degree blocks correspond to the leading forms of the polynomials
in $M$. Any finite set can be brought into canonical
form by applying Gaussian elimination and column permutations of terms
with same degree if necessary. In particular, the output of the $\lstabspan$
procedure can be easily brought into this form.  The following lemma
summarizes the basic properties of a set $M$ in canonical form:
\begin{lem}
\label{lem:stableIdeal}
Let $M$ be in canonical form and $L$-stabilized with $L=\T_{\leq d}^{n}$.
Let \(\T^{n}_{=d} \subseteq \vgen{M} + \T^{n}_{\leq d-1}\).
Then the following hold for all $i\in[d]$:
\begin{enumerate}
\item\label{item:1}
  \(\left. \igen{M}^{\leq i} \middle/ \igen{M}^{\leq i-1} \right. \cong
  \vgen{\LF(M_{i})}\)
\item \(\igen{M}^{\leq i} = \vgen{\bigcup_{j \leq i} M_{j}}\)
\item $\vgen{M_{i}}^{<i}=0$ and thus $\vgen{M_{i}}^{<i}\subseteq\vgen{\bigcup_{0\leq j\leq i-1}M_{j}}$
\end{enumerate}
\begin{proof}
We first show that $\vgen{M_{i}}^{<i}=0$ for all $i\in[d]$. Let
$i\in[d]$ be arbitrary and observe that each nonzero element $p\in M_{i}$
has degree $i$. As $M$ is in canonical form, the polynomials
in $M_{i}$ are interreduced (see the matrix in
Figure~\ref{fig:l-can-form} for Definition \ref{def:L-canonical-form})
and thus we also obtain each nonzero element $p\in\vgen{M_{i}}$ has
degree $i$.

By Corollary~\ref{cor:idealApprox}, \(\igen{M} \cap \vgen{L} = \vgen{M}\).
Hence \(\igen{M}^{\leq i} = \vgen{M}^{\leq i}\) for \(i \in [d]\).
Now the statements of the lemma are obvious consequences of \(M\)
being in canonical form.
\end{proof}
\end{lem}
The following lemma provides us a practical way to compute
the sizes of the degree components of degree-compatible order ideals,
which are the same for all order ideals of a given ideal.
\begin{lem}
\label{lem:stableRepresentation}
Let
$M$ be in canonical form and $L$-stabilized with $L=\T_{\leq d}^{n}$
and $d=\max_{m\in\partial\oi}\deg(m)$.
Let \(\T^{n}_{=d} \subseteq \vgen{M} + \T^{n}_{\leq d-1}\).
Further let $\oi$ be an order ideal of $\igen{M}$.
Then $\oi$ is degree-compatible if and only if \[
\lvert \oi^{=i} \rvert =\dim\vgen{L^{=i}}-\dim\vgen{\LF(M_{i})}\]
for every \(i \in [d]\).

\begin{proof}
In view of Definition \ref{def:indepSizes} it suffices to observe
that $\left. I^{\leq i} \middle/ I^{\leq i-1} \right.\cong\vgen{\LF(M_{i})}$
by Lemma~\ref{lem:stableIdeal} (\ref{item:1})
where \(I \coloneqq \igen{M}\).
\end{proof}
\end{lem}

We are ready to provide a reformulation of
the definition of order ideal polytopes,
which is better suited for actual computations, partly as
it only involves direct matrix operations
via replacing dimensions with ranks of subsets:
\begin{lem}
\label{lem:polytope-comp}

Let $M$ be $L$-stabilized and in canonical form with
$L=\T_{\leq d}^{n}$
and
\(\left. \vgen{M}^{\leq d} \middle/ \vgen{M}^{\leq d-1} \right. \cong
\vgen{\T^n_{= d}}\).
Then the \emph{order ideal polytope}
$P(M,L)\subseteq\cube L$
of \(\igen{M}\)
is given by the system of inequalities in \autoref{fig:OrderIdealPolytope}.

\begin{figure}
\begin{empheq}[box=\fbox]{align}
  z_{m_{1}} &\geq z_{m_{2}} \quad
  \forall m_{1},m_{2}\in L\colon m_{1}\mid m_{2}\label{eq:orderIdeal}\\
  \sum_{m\in L^{=i}}z_{m} &= \dim\vgen{L^{=i}}-\dim\vgen{\LF(M_{i})} \quad
  \forall i\in[d-1]\label{eq:degreeCompatibility}\\
  \sum_{m\in U}z_{m} &\geq
  \begin{multlined}[t]
    \lvert U \rvert -\rk{\tilde{U}}\label{eq:equiSteinitz}\\
    \forall i\in[d-1], U\subseteq L^{=i} \colon
    \lvert U \rvert=\dim\vgen{\LF(M_{i})}
  \end{multlined}
  \\
  z_{m} &\in [0,1] \quad \forall m\in L\nonumber
\end{empheq}

\caption[Order ideal polytope $P(M,L)$]{\label{fig:OrderIdealPolytope}%
  Order ideal polytope $P(M,L)$.
  In \autoref{eq:equiSteinitz}, the matrix
  $\tilde{U}$ is the induced sub-matrix of $\LF(M_{i})$ with column monomials
  only in $U$.}

\end{figure}

\begin{proof}
Let \(I \coloneqq \igen{M}\).
We successively transform the defining inequalities
of the order ideal polytope in \autoref{fig:OrderIdealPolytope2}
into the desired form of \autoref{fig:OrderIdealPolytope}.
The reformulation is mostly based on
\(\left. I^{\leq i} \middle/ I^{\leq i-1} \right.\cong\vgen{\LF(M_{i})}\)
from Lemma~\ref{lem:stableIdeal}(\ref{item:1}).

First, as \(\left. I^{\leq d} \middle/ I^{\leq d-1} \right. \cong
\vgen{\T^n_{= d}}\),
we can remove the variables \(z_m\) with \(m\) degree at least \(d\)
together with the inequalities involving them.
These variables are always zero.

Second, we replace all occurrences of \(\left. I^{\leq i} \middle/ I^{\leq i-1} \right.\) with
\(\vgen{\LF(M_{i})}\) (or simply \(\LF(M_{i})\)).
This almost results in the inequality system of
\autoref{fig:OrderIdealPolytope},
with the only difference that instead of \autoref{eq:equiSteinitz}
we have
\begin{multline}
  \sum_{m\in U}z_{m} \leq
    \dim\vgen{U'\cup\LF(M_{i})}-\dim\vgen{\LF(M_{i})}\label{eq:steinitz} \\
    \forall i\in[d-1],U\subseteq L^{=i}\colon \lvert U'
    \rvert=\dim\vgen{L^{=i}}-\dim\vgen{\LF(M_{i})},
\end{multline}
where we have deliberately replaced \(U\) with \(U'\).

We will show that the difference of \equationautorefname s
\ref{eq:degreeCompatibility} and \ref{eq:equiSteinitz}
is equal to \autoref{eq:steinitz} with the choice
$U'\coloneqq L^{=i}\setminus U$,
which has size $\lvert U' \rvert=\dim\vgen{L^{=i}}-\dim\vgen{\LF(M_{i})}$.
This will finish the proof.

Let $U\subseteq L^{=i}$ as above and compute the difference of \autoref{eq:degreeCompatibility} and \autoref{eq:equiSteinitz}.
We obtain \[
\sum_{m\in L^{=i}\setminus U}z_{m}\leq\dim\vgen{L^{=i}}-\dim\vgen{\LF(M_{i})}-
\lvert U \rvert + \rk{\tilde{U}}.\]
It is easy to see that $\lvert U \rvert
=\dim\frac{\vgen{L^{=i}}}{\vgen{L^{=i}\setminus U}}$.
We claim that it suffices to show that \[
\rk{\tilde{U}}=\dim\frac{\vgen{\LF(M_{i})\cup(L^{=i}\setminus U)}}{\vgen{L^{=i}\setminus U}}.\]
Indeed, using this we can rewrite the inequality as
\begin{equation*}
  \begin{split}
    \sum_{m\in L^{=i}\setminus U}z_{m} &\leq \dim\vgen{L^{=i}} -
    \dim\vgen{\LF(M_{i})} - \dim\frac{\vgen{L^{=i}}}{\vgen{L^{=i}\setminus U}}
    + \dim\frac{\vgen{\LF(M_{i}) \cup (L^{=i}\setminus U)}}{\vgen{L^{=i}
        \setminus U}}\\
    &= \dim\vgen{\LF(M_{i})\cup(L^{=i}\setminus U)}-\dim\vgen{\LF(M_{i})},
  \end{split}
\end{equation*}
which is \eqref{eq:steinitz} for $U'\coloneqq L^{=i}\setminus U$ as claimed.

We will show now that $\rk{\tilde{U}}=\dim\vgen{\LF(M_{i})\cup(L^{=i}\setminus U)}-\dim\vgen{L^{=i}\setminus U}$.
Let $B$ denote the matrix obtained when writing the elements in $\LF(M_{i})$
as rows and let $L^{=i}$ index the columns. Clearly, $\vgen{L^{=i}}=\bigoplus_{\ell\in L^{=i}}Ke_{\ell}$
as a vector space and $\vgen{U'}=\bigoplus_{\ell\in U'} Ke_{\ell}$
as a sub vector space with $U'\coloneqq L^{=i}\setminus U$. We obtain
$\left. \vgen{L^{=i}} \middle/ \vgen{U'} \right. \cong \bigoplus_{\ell\in
    L^{=i}\setminus U'} Ke_{\ell}$.
Note that $\vgen{\LF(M_{i})}\subseteq\vgen{L^{=i}}$ and thus we can
consider $\left. \vgen{\LF(M_{i})\cup U'} \middle/
  \vgen{U'} \right. \subseteq \left. \vgen{L^{=i}} \middle/ \vgen{U'} \right.$.
Now $\dim \left( \vgen{\LF(M_{i})\cup U'} \middle/ \vgen{U'} \right) =
\rk{\tilde{U}}$ where
$\tilde{U}$ is obtained from $B$ by removing the columns in $U'$.
We obtain
\begin{equation*}
  \begin{split}
    \rk{\tilde{U}} &= \dim\frac{\vgen{\LF(M_{i})\cup U'}}{\vgen{U'}}\\
    &= \dim\frac{\vgen{\LF(M_{i})\cup(L^{=i}\setminus
        U)}}{\vgen{L^{=i}\setminus U}}
  \end{split}
\end{equation*}
and thus the result follows.
\end{proof}
\end{lem}

\section{\label{sec:Computing-border-bases}Computing border bases using the
order ideal polytope}

In the following we explain how Theorem \ref{thm:polytope=order-ideal}
can be used to actually compute border bases for general degree-compatible
order ideals. We cannot expect to be able to compute a border basis
for \emph{any} degree-compatible order ideal, simply as such a basis
does not necessarily exist. Having the order ideal polytope at hand
something slightly more subtle can be done: By choosing a linear objective
function $c\in\Z^{\T^n}$ and optimizing it over the order ideal polytope,
we can actually \emph{search} for an order ideal with preferred monomials
in its support (see \autoref{sub:preferredBasis}). Having the
order ideal polytope available we can also count the number of degree-compatible
border bases that exist for a specific ideal.
Before we can address this application though, we will first show
how to obtain an $\oi$-border basis for $\oi \in \Lambda(I)$
where $I \subseteq \pring$ is a zero-dimensional ideal.

\subsection{\label{sub:ComputingBB}Computing border bases for $\oi \in \Lambda(I)$}

As the computation of the $L$-stable span of a set of generators
$M$ is independent of the actual chosen vector space basis (see Remark
\ref{rem:stableForVR}), we can adapt the classical border basis algorithm
(Algorithm \ref{alg:borderBasis}) to compute border bases for general
degree-compatible order ideal. We first determine the right computational
universe $L = \T^n_{\leq d}$ for some $d \in \N$ such that the associated
$L$-stable span $M$ contains all border bases. In a second step we optimize over
the order ideal polytope $P(M,L)$ and then perform the corresponding basis
transformation. We will first formulate the \emph{generalized border basis
  algorithm} by adding two steps after (\ref{enu:bbasisStopNew}) in Algorithm
\ref{alg:borderBasis} to the classical border basis algorithm, formulate the
missing parts, and then prove its correctness:

\begin{algorithm}{Generalized border basis algorithm}{\bbasis}%
{$F$ a finite generating set of a zero-dimensional ideal.}%
{$\mathcal{G}$ a border basis of the ideal.}
\label{alg:borderBasisNew}

\item Let $d\coloneqq \max_{f\in F}\{\deg(f)\}$ and put $L\coloneqq \T_{\leq d}^{n}$.
\item \label{enu:bbasisStep2_1}$V=\{v_{1},\dots,v_{r}\}\coloneqq \lstabspan(F,L)$.
\item \label{enu:bbasisStopNew}
If \(\T^{n}_{= d} \nsubseteq \LT(V)\)
then set $d\coloneqq d+1$
and put $L\coloneqq \T_{\leq d}^{n}$ and go to step (\ref{enu:bbasisStep2_1}).
\item \label{enu:calcOrderIdeal} Choose $\oi \in \Lambda(\igen{V})$ ($\Leftrightarrow z \in P(V,L) \cap \Z^L$ and $\oi = \oi(z)$).  
\item \label{enu:basisTransformation} Let $\mathcal{G} \coloneqq   \basisT(V,\oi)$.
\end{algorithm}

Note that step (\ref{enu:bbasisStopNew})
is a convenient way to quickly check whether $V$ is already
$L$-stabilized and if it contains all degree-compatible order ideals. In this
augmented algorithm we added the steps (\ref{enu:calcOrderIdeal}) and
(\ref{enu:basisTransformation}). The first step will be extensively discussed
in \autoref{sub:preferredBasis} as there are various ways to determine $\oi
\in \Lambda(\igen{V})$ and this is precisely one of the main features, i.e.,
to choose the order ideal more freely. Note that by
Theorem~\ref{thm:polytope=order-ideal} we
already know that $\oi$ does support an $\oi$-border basis of $\igen{F}$ (as
$\igen{F} = \igen{V}$) and our task is now to actually extract this basis
from $V$. This extraction is performed in step
(\ref{enu:basisTransformation}).

Let $A$ be a matrix representing a set of polynomials $M$ where the columns
correspond to the monomials in some fixed ordering.
Let the \emph{head} (short: $\head(a)$) of a row $a$ of $A$ be
the left-most monomial in
the matrix representation whose coefficient is non-zero. Note that the notion
of \emph{head} replaces the notion of \emph{leading term} of a polynomial as we
do not (necessarily) have a term ordering anymore.  The main idea is to
reorder the columns of $V$ and then to
bring $V$ into a reduced row echelon form such that
no $m \in \oi$ is head of a row of the resulting matrix
\textemdash\ a classical basis transformation:

\begin{lem}
  \label{lem:basisT}
  Let $L = \T^n_{\leq \ell}$ with $\ell\in \N$,
  let $V$ be a finite set of polynomials satisfying
  \(\vgen{V} = \igen{V} \cap \vgen{L}\)
  and
  let $\oi = \{t_1, \dots, t_\mu\}$ be an order ideal with
  $\partial \oi \subseteq L$
  and $\oi \in \Lambda(\igen{V})$.
  Then Algorithm~\ref{alg:basisT} returns
  the $\oi$-border basis $\mathcal G$ of $\igen{V}$.
\begin{proof}
First, the algorithm finds \(\ell\) from \(M\).
As $\oi \in \Lambda(\igen{V})$ we have that $\oi$ supports a border
basis and in particular we have $\pring = \igen{V} \oplus \vgen{\oi}$,
and hence \(\vgen{L} = \vgen{V} \oplus \vgen{\oi}\)
by the modular law.
We will now show
that \hyperref[enu:bbform]{Condition (\ref*{enu:bbform})} of
Definition~\ref{def:bbasis} is
satisfied. This in turn follows from the fact that
the algorithm creates every element $g_j$ of $\mathcal G$ to have the form
\[g_j =b_{j}-\sum_{i=1}^{\mu}\alpha_{ij}t_{i}\] with $\alpha_{ij} \in K$. By
construction $b_j \in \partial \oi$ and thus $\mathcal G$ is an $\oi$-border
basis of $\igen{V}$.
\end{proof}
\end{lem}

\begin{algorithm}{Basis transformation algorithm}{\basisT}%
{$V,\oi$ as in Lemma \ref{lem:basisT}.}%
{$\mathcal G$ as in Lemma \ref{lem:basisT}.}
\label{alg:basisT}

\item Set \(\ell \coloneqq \max_{m \in M} \deg(m)\).
\item
  Permute the columns of $V$ such that
  $t$ is right of $m$ in the matrix representation of $V$
  for all $m \in \T^n_{\leq \ell}$ and $t \in \oi$.
\item
  Reduce $V$ using Gaussian elmination: like Algorithm~\ref{alg:gaussEl},
  but use \(\head\) instead of \(\LT\) and the coefficient of the head
  instead of \(\LC\).
  Let $\mathcal G'$ be the result.
\item Let $\mathcal G \coloneqq \{g \in \mathcal G' : \head(g) \in \partial\oi \}$.
\item Return $\mathcal G$.
\end{algorithm}
Note that Algorithm~\ref{alg:basisT} works for \emph{any} order ideal
that supports a border basis of \(\igen{V}\), i.e., also those that
are not necessarily degree-compatible. When the order ideal is known to be degree-compatible,
it is enough to do the permutations in each degree block
in the first step,
and then use Algorithm~\ref{alg:gaussEl} in the second step.

We will show now that Algorithm~\ref{alg:borderBasisNew} computes an $\oi$-border basis for $\oi \in \Lambda(I)$.

\begin{prop}
\label{pro:bbasisNew}Let $F=\{f_{1},\dots,f_{s}\}\subseteq \pring$
be a finite set of polynomials that generates a zero-dimensional ideal
$I=\igen{F}$. Then Algorithm \ref{alg:borderBasisNew} computes
the $\oi$-border basis $\mathcal{G}$ of $I$ for any (chosen) $\oi \in \Lambda(I)$.
\begin{proof}
Whenever we reach step (\ref{enu:calcOrderIdeal}) in
Algorithm~\ref{alg:borderBasisNew}, we have that $V$ is $L$-stabilized for
some $L = \T^n_{\leq d}$ with $d\in \N$ and it contains all degree-compatible
order ideals supporting a border basis, i.e., all $\oi \in \Lambda(I)$.
Observe that $I = \igen{F} = \igen{V}$ and thus, by
Lemma~\ref{lem:basisT}, it follows that $\mathcal G$ is indeed an
$\oi$-border basis of $\igen{F}$.
Note that step~(\ref{enu:bbasisStopNew}) ensures
$\partial \oi \subseteq L$
via \(\left. \vgen{V}^{\leq d} / \vgen{V}^{\leq d-1} \right. \cong
 \vgen{\T^{n}_{= d}} = \vgen{L^{= d}}\).
\end{proof}
\end{prop}

An improved version of the border basis algorithm has been also considered in
\cite{kehrein2006cbb}. Basically, the improvement can be traced back to
considering more restricted computational universes $L$ that arise from
choosing $L$ to be the smallest order ideal that contains the support of the
initial system and then successively extending it using the $+$
operation. This improvement due to the restriction of the computational universe
$L$ cannot work in our setting anymore: Suppose that $V$ is $L$-stabilized with
respect to some computational universe $L \neq \T^n_{\leq d}$ for all $d
\in \N$ (i.e., \(L\) is not obtained by bounding the total degree of
the monomials in \(\T\)) and contains the order ideal that is induced by the chosen
degree-compatible term ordering in the classical border basis algorithm. Then
the associated polytope $P(V,L)$ would only contain a \emph{subset} of all
possible degree-compatible order ideals, as we might be lacking monomials that
we need to represent certain alternative choices of $\oi$. If a subset of all admissible 
degree-compatible order ideals is sufficient, then the same optimizations can be applied though.

\subsection{\label{sub:preferredBasis}Computing preferred border bases}
Let $V$ and $L$ be as obtained after step (\ref{enu:bbasisStopNew}) in
Algorithm \ref{alg:borderBasisNew}.  As shown in
Theorem~\ref{thm:polytope=order-ideal} and
Lemma~\ref{lem:polytope-comp}, the order ideal polytope $P(V,L)$ characterizes
all degree-compatible order ideals that support a border basis of
$\igen{V}$. Every $z \in P(V,L) \cap \Z^L$ induces an order ideal $\oi(z)$
which supports an $\oi(z)$-border basis of $\igen{V}$. This
also shows that we cannot expect that every order ideal $\oi$ supports an
$\oi$-border basis of $\igen{V}$ as the characterization is one-to-one. 
The natural question is therefore how to
\emph{specify} which order ideal should be computed, i.e., which $z \in P(V,L)
\cap \Z^L$ to choose. As we cannot always get what we would like to have, it
suggests itself to specify a \emph{preference}, i.e., which monomials we would like
to be contained in $\oi(z)$ and which ones we would rather not.
As the coordinates
of $z$ are in direct correspondence with the monomials in $L$ we can define:

\begin{defn}
  A \emph{preference} is a vector $c \in \Z^L$ which assigns a weight to each
  monomial $m \in L$. If $z\in P(V,L) \cap \Z^L$, then $cz$ is the
  \emph{score} or \emph{weight} of $z$.
\end{defn}

As $P(V,L) \subseteq \cube{n}$ is a polytope we can optimize over $P(V,L) \cap
\Z^L$ and compute an element $z_0 \in P(V,L) \cap \Z^L$ that has maximal
score, i.e., we can compute $z_0 \in P(V,L) \cap \Z^L$ such that
\[cz_0 = \max\{ cz \mid z \in P(V,L) \cap \Z^L \}.\] In this sense a
preference is an indirect way of specifying an order ideal. For
certain choices of $c \in \Z^L$ though it can be hard to compute an
order ideal that maximizes the score as we will show now.

\section{Complexity of finding preferred order ideals}
\label{sec:complexity-finding-ord-ideal}

In this section,
we show that finding a weight optimal, border basis supporting order ideal of
a zero-dimensional ideal given by generators is \NP-hard
(Theorem~\ref{thm:MaxOrderIdeal-NPhard}).
Note that this also translates to large
classes of other choice functions as $P(V,L) \cap \Z^L$ is a 0/1 polytope and
thus its extremal points are given by an inequality description. So any choice
function that implicitly asks for a $k$-clique (which we will use in our reductions) can be replaced by the
corresponding linear function and hardness also follows in this case.
The hardness for computing a weight optimal order ideal is
unexpected in the sense that we merely ask for a basis transformation. On the
other hand it highlights the crucial role of order ideals in describing the
combinatorial structure of the ideal.

As an immediate consequence it follows that it is rather unlikely that we can
obtain a good characterization of the integral hull $\conv{P(V,L) \cap \Z^L}$
and we will not be able to compute degree-compatible order ideals that support
a border basis and have maximum score efficiently unless $\Pclass=\NP $. This shows
that not only computing the necessary liftings of the initial set of
polynomials via the $\lstabspan$ procedure is hard but also actually
determining an optimal choice of an order ideal once an $L$-stable span has
been computed. As mentioned before, this is in some sense surprising as the actual interference has
been already performed at that point and we are only concerned with choosing a
\emph{nice} basis. From a practical point of view this is not too problematic
as, although \NP-hard, computing a weight optimal order ideal is no harder
than actually computing the $\lstabspan$ in general. For bounds on the degree
$d \in \N$ needed to compute border bases see, e.g., \cite[Lemma
2.4]{deLoera2009ecps}; the border basis algorithm generates the
Nullstellensatz certificates and is therefore subject to the same
bounds. Further, state-of-the-art mixed integer programming solvers that can
solve the optimization problem such as \textprogram{scip}
(\cite{achterberg2009scip}), \textprogram{cplex} (\cite{cplex200811}), or
\textprogram{gurobi} (\cite{gurobi}) can handle instance sizes far beyond the
point for which the actual border bases can be computed. Very good solutions
can also be generated using simple local search schemes starting from a
feasible order ideal derived from a degree-compatible term ordering.

\subsection{Fast without constraint}
\label{sec:fast-with-constr}

\begin{figure}
  \centering
\includegraphics[scale=0.55]{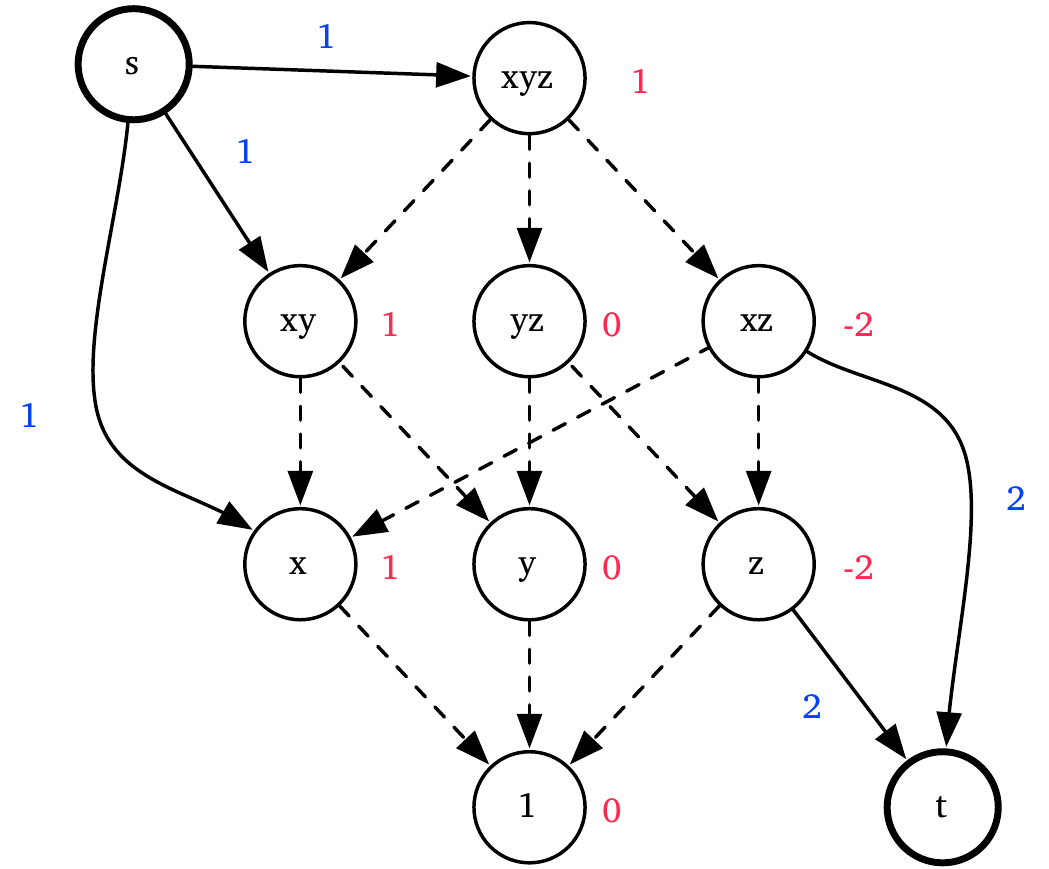}
\caption[Order ideal computation as minimum cut problem]%
{\label{fig:orderIdeal}%
Order ideal computation as minimum cut problem.
Capacities denoted next to the arc and weight $c_u$ denoted next to the
respective node $u$. The dashed arcs have capacity $\infty$.}
\end{figure}

Determining an order ideal of maximum
score in a computational universe \(L\) without having any additional
constraints can be done in time
polynomial in $\lvert L \rvert $, as was shown in
\cite{picard1976maximal}.
One simply transforms it into a \emph{minimum cut problem} in
graphs (see, e.g., \cite{schrijver1986theory}): Let $c \in \Z^L$ be a
preference vector. Define a
graph $\Gamma \coloneqq (V,E)$ with $V \coloneqq L \cup \{s,t\}$ and $\tilde E \coloneqq  \{ (u,v)
\mid u,v \in L\text{ and } v \mid u \}$, i.e., whenever $v \mid u$ we add an arc from
$u$ to $v$.
In fact, it is enough to have an arc when \(u = v x\) for some variable \(x\), i.e., to consider the transitive reduction of $\tilde E$.
Define \[E \coloneqq  \tilde E \cup \{ (s,u) \mid u \in L, c_u > 0\} \cup
\{ (u,t) \mid u \in L, c_u < 0 \}.\] Further we set all the capacities of the
arcs with both vertices in $L$ to $\infty$, for any arc $(s,u)$ with $u \in L$
we set the capacity to $c_u$, and for any arc $(u,t)$ with $u \in L$ we set
the capacity to $\lvert c_u \rvert $; let $\kappa(u,v)$ denote the capacity of the arc
$(u,v) \in E$. An example is depicted in \autoref{fig:orderIdeal}.

For $U,W \subseteq V$, we define $C(U,W) \coloneqq  \sum_{(u,w) \in U \times W} \kappa(u,w)$. 
A \emph{(s,t)-cut} $(S,\bar S)$ is a partition $S \stackrel{\cdot}{\cup} \bar S=V$ of the vertices of $V$ with $s \in S$ and $t \in \bar S$ and the \emph{weight} of the cut is $C(S,\bar S)$. 

We would like to compute an order ideal contained in $L$ with maximum score:
\[ \max \left\{ \sum_{u \in \oi} c_u \middle| \oi \subseteq L \text{ order ideal} \right\}.\] 
Observe that, if $(S, \bar S)$ is a cut in $\Gamma$ of finite weight,
there exist no arc $(u,v) \in \tilde E$ with $u \in S$ and $v \in \bar S$. Therefore, $(S, \bar S)$ is a cut in $\Gamma$ of finite weight if and only if $S \setminus \{s\}$ is an order ideal. We can therefore rewrite the optimization problem as follows:
\begin{equation*}
  \begin{split}
    \max \left\{ \sum_{u \in \oi} c_u \middle| \oi \subseteq L \text{ order
        ideal} \right\} &= \max \{ C(\{s\},\oi ) - C(\oi, \{t\}) \mid
    \oi \subseteq L \text{ order ideal}\} \\
    &= \max \{ C(\{s\}, L ) - C(\{s\},L \setminus \oi ) - C(\oi, \{t\}) \mid
    \oi \subseteq L \text{ order ideal} \} \\
    &= C(\{s\}, L ) - \min \{ C(\{s\},L \setminus \oi ) + C(\oi, \{t\}) \mid
    \oi \subseteq L \text{ order ideal} \}  \\
    &= C(\{s\}, L ) - \min \{ C(\{s\} \cup \oi , (L \setminus \oi) \cup \{t\}
    ) \mid \oi \subseteq L \}.
  \end{split}
\end{equation*}
The last line asks for a minimum weight cut in the graph $\Gamma$. Note that we can indeed drop the condition that $\oi$ has to be an order ideal as it is guaranteed implicitly by all finite weight cuts as explained above. It is well-known that the computation of a minimum cut can be performed in time polynomial in the number of vertices and arcs (see \cite{wolsey1999integer}). Thus we can indeed compute an order ideal $\oi$ of maximum score efficiently in this case.

\subsection{\NP-hard with constraints}
\label{sec:np-hard-constraints}

So far we did not include the additional requirements as specified by the
order ideal polytope (see \autoref{fig:OrderIdealPolytope}), in order to
obtain degree-compatible order ideals that do actually support a border basis
of the ideal $I$ under consideration. Unfortunately, when including these
additional requirements, the problem of computing an order ideal of maximum
score is \NP-hard as we will show in the following.

We will show \NP-hardness by a reduction from the \nameref{pr:MaxClique} problem, which is well known to be \NP-complete
(see, e.g.,
\cite{garey1979cai} or \cite[GT22]{crescenzi-compendium}). Given a graph $\Gamma = (V,E)$, recall that a
\emph{clique} $C$ is a subset of $V$ such that for all distinct $u,v \in C$ we have
$(u,v) \in E$.

\begin{problem}{Max-Clique}
  \label{pr:MaxClique}
  Let $\Gamma = (V,E)$ be a graph. Determine the maximum size of a clique $C$
  contained in $\Gamma$.
\end{problem}

We will in particular use the following variant:

\begin{problem}{$k$-Clique}
  \label{pr:kClique}
  Let $\Gamma = (V,E)$ be a graph. Determine whether $\Gamma$ contains a
  clique $C$ of size $k$.
\end{problem}

Note that if $\Gamma$ contains a clique of size $k$ for some $k\in \N$, then
so does it for any $1 \leq l \leq k$. Further the maximum size of a clique is
bounded by $\lvert V \rvert $. It follows that already testing whether $\Gamma$ contains a
clique $C$ of size $k$ for some given $k \in \N$ has to be \NP-complete, as
otherwise, we could solve the \nameref{pr:MaxClique} problem with
$O(\log_2 \lvert V \rvert )$
calls of the algorithm that solves the test variant via binary search. In
\cite[Discussion after Definition 3.2]{hochbaum2000performance} it was
indicated that determining a maximum weight order ideal of a pre-defined size
is \NP-hard by a reduction from \nameref{pr:MaxClique}. While already indicating
hardness due to one additional cardinality constraint, this is a slightly
different problem than the one that we are facing here:
in addition to the size constraint, the order ideal has to be degree-compatible.
Thus we have a
cardinality/degree-constraint \emph{for each degree} of the monomials
(Constraints \eqref{eq:degreeCompatibility} and \eqref{eq:steinitz} in
\autoref{fig:OrderIdealPolytope}). Further, the degree constraints are
not completely independent of each other. So it could be \emph{a
  priori} perfectly possible that this particular restriction of the problem
can be actually solved efficiently, which is not the case as we will see soon.

We consider the following optimization problem: 

\begin{problem}{Max-Bounded order ideal}
  \label{pr:MaxBoundedOrderIdeal}
  Let $L = \T^n_{\leq d}$ for some $d \in \N$ and $c \in \Z^L$ be a
  preference. Further let $d_i \in \N$ for $i \in [d-1]$. Determine an order
  ideal $\oi_0 \subseteq L$ with $\lvert \oi_0^{=i} \rvert =d_i$ for all $i
  \in [d-1]$ and maximum score, i.e., \[c(\xi^{-1}(\oi_0)) =
  \max\{c(\xi^{-1}(\oi)) \mid \oi \subseteq L \text{ order ideal and } \lvert
  \oi^{=i} \rvert =d_i \text{ for all } i \in [d-1] \}.\]
\end{problem}

By a reduction from \nameref{pr:kClique} we obtain:

\begin{thm}
\label{thm:maxbboi}
\nameref{pr:MaxBoundedOrderIdeal} is \NP-hard.
\begin{proof}
  Let $\Gamma = (V,E)$ be an arbitrary simple graph and choose $k \in [\lvert
  V \rvert]$. We consider the polynomial ring $K[x_v \mid v \in
  V]$. Let $L = \T^{\lvert V \rvert }_{\leq 2}$. Further choose $d_1 =
  k$ and $d_2 = k(k+1)/2$. Define
\[c_m = 
\begin{cases}
1, \text{ if } m = x_u x_v \text{ and } (u,v) \in E; \\
0, \text{ otherwise}. 
\end{cases}
\]
A vertex $v \in V$ is represented by the degree-one monomial $x_v$ and each edge
$(u,v) \in E$ is represented by the degree-two monomial $x_ux_v$. In some sense
we extended the graph $(V,E)$ to the complete graph $K_{\lvert V \rvert }$
that has $\lvert V \rvert (\lvert V \rvert -1)/2$ edges as we consider all
possible degree-two monomials that correspond to the edges; those edges that
do not belong to $E$ have weight $0$ though.

The order ideals satisfying the size constraints are
those consisting of the monomials \(1\), \(x_{v_i}\),
\(x_{v_i} x_{v_j}\)
for \(1 \leq i, j \leq k\) for some distinct vertices
\(v_1\),\dots,\(v_k\).
The score of the order ideal is the number of edges in the subgraph
spanned by the \(x_{v_i}\).

Thus when maximizing $c$ we ask for a
\(k\)-vertex subgraph with maximal number of edges.
The subgraph is given by the degree-one monomials contained in the
order ideal, which we denote by $\oi$, i.e.,
\[C \coloneqq \{v \in V \mid x_v \in \oi \}.\]

So the maximum score is equal to $k(k-1)/2$ if and only if $\Gamma$ contains a clique of size $k$.
Thus we can solve
\nameref{pr:kClique} efficiently if we can solve
\nameref{pr:MaxBoundedOrderIdeal} efficiently and so
the latter has to be \NP-hard.
\end{proof} 
\end{thm}

Finally, it remains to show that
for every graph $\Gamma = (V,E)$ and $k\in [\lvert V
\rvert]$ there exists a
system of polynomials
$F_{\lvert V \rvert,k} \subseteq K[x_v \mid v \in V]$
spanning a zero-dimensional ideal
such that solving the \nameref{pr:MaxOrderIdealOfIdeal} problem for
$F_{\lvert V \rvert,k}$ solves the \nameref{pr:kClique} problem for $\Gamma$.
For this, we construct an ideal encoding all $k$-cliques of the
complete graph on $n$ vertices: Let $n \in \N$ and $k \in [n]$ and define 
\begin{gather*}
  v_j \coloneqq \sum_{i \in [n]} i^j x_i, \\
  F_{n,k} \coloneqq \{ v_j \mid j \in [n-k] \} \cup \T^n_{=3}.
\end{gather*}
We consider the ideal generated by \(F_{n,k}\).
We show that its order ideals are in one-to-one correspondence with the
\(k\)-element subsets of the set of \(n\) variables \(x_1,\dots,x_n\)
as stated in the following
lemma.

\begin{lem}
\label{thm:thePolySystem}
Let $n \in \N$ and $k\in [n]$. Then $F_{n,k}$ generates a zero-dimensional
ideal such that $\oi \in \Lambda \left( \igen{F_{n,k}} \right)$ if and only if
$\oi^{=1} \subseteq \T^n_{=1}$ with $\lvert \oi^{=1} \rvert = k$, $\oi^{=2} = \{xy
\mid x,y \in \oi^{=1}\}$, and $\oi^{=\ell} = \emptyset$ for all $\ell \geq 3$.
\begin{proof}
We will first characterize the vector space $\left. K[x_1, \dots, x_n]
  \middle/ \igen{F_{n,k}} \right.$. Observe that the coefficient matrix $A
\coloneqq (v_j)_{j \in [n-k]}$ is actually a Vandermonde matrix and in
particular every square submatrix of $A$ is invertible. Furthermore, the
polynomials $v_j$ with $j \in [n-k]$ are homogeneous of degree one. Thus, for
any $k$ variables of $\{x_1, \dots, x_n\}$, without loss of generality say,
$x_1, \dots, x_k$, we have that $\{x_1, \dots, x_k, v_1, \dots , v_{n-k}\}$ is
a basis for the homogeneous polynomials of degree one. This is just
another way of saying that removing the columns belonging to $x_1, \dots, x_k$
from $A$, the resulting square submatrix is invertible.  It follows
\begin{equation*}
  \begin{split}
    \vgen{\oi} \cong \left. K[x_1, \dots, x_n] \middle/ \igen{F_{n,k}} \right.
    &\cong \left. K[x_1, \dots, x_k, v_1, \dots , v_{n-k}] \middle/
      \igen{F_{n,k}} \right. \\
    &\cong \left. K[x_1, \dots, x_k] \middle/ \igen{\T^k_{=3}} \right. .
  \end{split}
\end{equation*}
As the substitution preserves degrees, homogeneity, etc., it follows, that any
degree-compatible order ideal has to have $\lvert \oi^{=1} \rvert = k$,
$\oi^{=2} = \{xy \mid x,y \in \oi^{=1}\}$, and $\oi^{=\ell} = \emptyset$ for
all $\ell \geq 3$.

The other direction follows immediately as each order ideal $\oi$ with $\lvert
\oi^{=1} \rvert = k$, $\oi^{=2} = \{xy \mid x,y \in \oi^{=1}\}$, and
$\oi^{=\ell} = \emptyset$ for all $\ell \geq 3$ is actually a
degree-compatible order ideal such that $K[x_1,\dots, x_n] = 
\igen{F_{n,k}} \oplus \vgen{\oi}$ as vector spaces, by the argumentation above.
\end{proof}
\end{lem}

Note that the order ideals of $F_{n,k}$ indeed correspond to the
$k$-cliques of the complete graph on $n$ vertices: If $\oi \in
\Lambda(F_{n,k})$, then $\oi^{=1} = \{x_{i_1}, \dots, x_{i_k}\}$ and
$x_{i_j}x_{i_l} \in \oi^{=2}$ if and only if $x_{i_j},x_{i_l} \in \oi^{=1}$. If we now
remove all elements of the form $x_{i_j}^2$ with $x_{i_j} \in \oi^{=1}$, and
there are $k$ of those, then \[ \lvert \oi^{=2} \setminus \{x_{i_j}^2 \mid
x_{i_j} \in \oi^{=1}\} \rvert = \frac{k(k-1)}{2},\]
the size of a $k$-clique. We are ready to state the main result of
\hyperref[sec:Computing-border-bases]{this section}.
Consider the following problem: 

\begin{problem}{Max-Order ideal of ideal}
  \label{pr:MaxOrderIdealOfIdeal}
  Let $M \subseteq \pring$ be a system of polynomials generating a
  zero-dimensional ideal and let $c\in \Z^{\T_n}$ be a preference on the monomials. Compute
  an order ideal $\oi$ supporting an $\oi$-border basis of $\igen{M}$ with
  maximum score with respect to $c$.
\end{problem}

\begin{thm}
  \label{thm:MaxOrderIdeal-NPhard}
  \nameref{pr:MaxOrderIdealOfIdeal} is \NP-hard.
\begin{proof}
The proof is by a reduction from the \NP-hard \nameref{pr:kClique} along the
lines of the
proof of Theorem~\ref{thm:maxbboi}.  Let $\Gamma = (V,E)$ be a graph with $n
\coloneqq \lvert V \rvert$ and $k \in [n]$ be an instance of
\nameref{pr:kClique}.  We consider $M \coloneqq F_{n,k}$ and
define $c \in \Z^{\T_{\leq
    3}^n}$ via
\[c_m = 
\begin{cases}
1, &\text{if } m = x_u x_v \text{ and } (u,v) \in E; \\
0, &\text{otherwise}. 
\end{cases}
\]
By Lemma~\ref{thm:thePolySystem}, we have that the degree-compatible order
ideals of $\igen{M}$ are in one-to-one correspondence with the $k$-cliques of
the complete graph on $n$ vertices.  Similarly to the proof of
Theorem~\ref{thm:maxbboi},
the score of an order ideal is just the number of edges between the
corresponding \(k\) vertices of the graph,
so the maximum score is \(k(k-1)/2\) if and only if the graph contains a
\(k\)-clique.
Thus, we obtain a test whether $\Gamma$
contains a clique of size $k$.
\end{proof}
\end{thm}

\section{\label{sec:Computational-results}Computational results}
We performed a few computational tests to verify the practical feasibility of
our method. All computations were performed with \textprogram{CoCoA 4.7.5}
(\cite{team2009cocoa}) and \textprogram{scip 1.1.0} (\cite{achterberg2009scip}) on
a 2 GHz Dual Core Intel machine with 2 GB of main memory. We performed
computations on various sets of systems of polynomial equations. The employed
methodology was as follows. We first computed a border basis using the
classical border basis algorithm. From the last run of the algorithm we
extracted the $L$-stabilized span and brought it into canonical form. We
generated the constraints \eqref{eq:degreeCompatibility} from the order ideal
that we obtained; from the $L$-stabilized span in matrix from, we generated
the constraints \eqref{eq:equiSteinitz}. As we needed to get access to the
$L$-stable span computed in the last round of the border basis algorithm, we
had to use an implementation of the border basis algorithm in
\textprogram{CoCoA-L} which is
slower in terms of speed compared to a C or C++ implementation. We then
transcripted these constraints into the CPLEX LP format which served as input
for \textprogram{scip}. For the optimization we chose various preference
functions.  One was a random function, and the other one was chosen with the
intent to make the optimization particularly hard by giving monomials deep in
the order ideal negative weights and assigning positive weights for the outer
elements.

In all cases the optimization, i.e., the computation of the weight optimal
order ideal was performed in less than a second, whereas the actual
calculation of the initial border bases was significantly more time consuming.
As indicated before, this is not unexpected as the computation of the
$L$-stable span is significantly more involved than computing a weight optimal
order ideal (in the worst case double exponential vs.~single exponential).
When computationally feasible we also counted all feasible order ideals with
\textprogram{scip}, which basically means enumerating all feasible solutions.
This in fact is equivalent to optimizing \emph{all} potential preference
functions simultaneously and thus emphasizes the
computational feasibility.

We considered $7$ systems of polynomials of various complexity in terms of
number of variables and order ideal degree. We would have liked to test
significantly larger instances but we were not able to compute the initial
border basis (or more precisely the $L$-stable span), neither with our
implementation nor with the C/C++ implementation of \textprogram{CoCoA 5}
(\textalgorithm{BBasis5}). This shows once more that the limiting factor is the actual computation of the $L$-stable span.

In Table~\ref{tab:results} we report our results.
\begin{table}[htdp]
\begin{center}
\begin{tabular}{|>{\centering}p{4cm}|c|c|c|c|}
\hline
 polynomial system & order ideal signature & optimization [s] & counting [s] & \# order ideals \\
\hhline{|=|=|=|=|=|}
$x^3, xy^2+y^3$ & $(1,3,1,1,1)$ & < 0.01 & 0.02 & 3 \\
\hline
vanishing ideal of the points \((0,0,0,1)\), \((1,0,0,2)\), \((3,0,0,2)\),
\((5,0,0,3)\), \((-1,0,0,4)\), \((4,4,4,5)\), \((0,0,7,6))\).
 & $(1,4,2) $ & < 0.01 & 0.02 & 45 \\
\hline
$x+y+z-u-v$, $x^2-x$, $y^2-y$, $z^2-z$, $u^2-u$, $v^2-v$ & $(1,4,5) $ & < 0.01 & 0.35 & 1,260 \\
\hline
$x+y+z-u-v$, $x^3-x$, $y^3-y$, $z^2-z$, $u^2-u$, $v^2-v$ & $(1,4,7,6)$ & 0.02 & 51.50 & 106,820 \\ 
\hline
$x+y+z-u-v$, $x^3-x$, $y^3-y$, $z^3-z$, $u^2-u$, $v^2-v$ & $(1,4,8,9)$ & 0.02 & 53.00 & 108,900 \\
\hline
$x+y+z-u-v$, $x^3-x$, $y^3-y$, $z^3-z$, $u^3-u$, $v^2-v$ & $(1,4,9,12,9)$ & 0.08 & 300.00* & > 1,349,154 \\
\hline
$x+y+z-u-v+a$, $x^2-x$, $y^2-y$, $z^2-z$, $u^2-u$, $v^2-v$, $a^2-a$ & $(1,5,9) $ & < 0.01 & 8.68 & 30,030 \\
\hline
\end{tabular}
\medskip
\caption{\label{tab:results}Computational results.  The first column contains
  the considered polynomial system. The second column contains the degree
  vector of the order ideal, i.e.,
  ${(\dim \left. I^{\leq i} \middle/ I^{\leq i-1} \right.)}_i$ starting
  with $i=1$ and \(I^{\leq0} \coloneqq 0\). The third column contains the average
  time (in seconds) needed to optimize a random preference over the order
  ideal polytope (we performed 20 runs for each system). The fourth column
  contains the time (in seconds) needed to count all admissible
  degree-compatible order ideals and the last column contains the actual
  number of admissible degree-compatible order ideals. The \lq{}*\rq{}
  indicates that the counting had been stopped after 300 seconds. The number
  of order ideals reported in this case is the number that have been counted
  up to that point in time.}
\end{center}
\end{table}

\section{\label{sec:Concluding-remarks}Concluding remarks}
We provided a way to characterize all degree-compatible order ideals that
support a border basis for a given zero-dimensional ideal $I$ by borrowing
from combinatorial optimization and in particular polyhedral theory. We
established a one-to-one correspondence of the integral points of a certain
polytope, the order ideal polytope, and those degree-compatible order ideals
that support a border basis of $I$. This connection in particular links the ideals to their combinatorial structure of the factor spaces. 

Using this polytope we adapted the
classical border basis algorithm in order to be able to compute border basis
for general degree-compatible order ideal based on a preference ordering on
terms contained therein. Effectively, the algorithm can be used for any integral point contained in the order ideal polytope. 

We also showed that computing a border basis for a preference on the monomials one might want to have included in the order ideal is \NP-hard and thus we cannot expect to be able to efficiently compute preferred order ideals in general although it is merely  a basis transformation. On the other hand, this is not restricting the applicability of our method in any practical application because the preceding computation of the $L$-stable span dominates in terms of computational complexity. We finally presented a few computational results showing the applicability of our method for actual computations. 

\bibliographystyle{plainurl}
\bibliography{gcNotesBib0}

\end{document}